\documentclass[12pt]{article}
\usepackage{amsmath}
\usepackage{graphicx}
\usepackage{CJK}
\usepackage{multirow}
\usepackage{makecell}

\usepackage{subfigure}
\usepackage[justification=centering]{caption}
\input{amssym.tex}

\def\bc{\begin{center}}       \def\ec{\end{center}}
\def\ba{\begin{array}}        \def\ea{\end{array}}
\def\be{\begin{equation}}     \def\ee{\end{equation}}
\def\bea{\begin{eqnarray}}    \def\eea{\end{eqnarray}}
\def\beaa{\begin{eqnarray*}}  \def\eeaa{\end{eqnarray*}}

\def\mathbb{\Bbb}

\begin{document}
\baselineskip 18pt
\centerline {\bf \large On the number of limit cycles for polycycles}
\vskip 0.1 true cm
\centerline {\bf \large  $S^{(2)}$ and $S^{(3)}$ in quadratic Hamilton systems }
\vskip 0.1 true cm
\centerline {\bf \large under perturbations of piecewise smooth polynomials}

\vskip 0.3 true cm

\centerline{\bf  Jiaxin Wang, Liqin Zhao$^{*}$}
 \centerline{ School of Mathematical Sciences, Beijing Normal University,} \centerline{Laboratory of Mathematics and Complex Systems, Ministry of
Education,} \centerline{Beijing 100875, The People's Republic of China}

\footnotetext[1]{
This work was supported by NSFC(11671040)   \\ * Corresponding author.
E-mail:  zhaoliqin@bnu.edu.cn (L. Zhao).}
\vskip 0.2 true cm
\noindent{\bf Abstract} In this paper, by using Picard-Fuchs equations and Chebyshev criterion, we study the bifurcate of limit cycles for quadratic Hamilton system $S^{(2)}$ and $S^{(3)}$: $\dot{x}= y+2axy+by^2$, $\dot{y}=-x+x^2-ay^2$ with $a\in(-\frac{1}{2},1)$, $b=(1-a)(1+2a)^{1/2}$ and $a=1$, $b=0$ respectively, under perturbations of piecewise smooth polynomials with degree $n$. The discontinuity is the line $y=0$. We bound the number of zeros of first order Melnikov function which controls the number of limit cycles bifurcating from the center. It is proved that the upper bounds of the number of limit cycles for $S^{(2)}$ and $S^{(3)}$ are respectively $25n+161$ $(n\geq3)$ and $24n+126$ $(n\geq3)$ (taking into account the multiplicity).

\noindent{\bf Keywords} limit cycle; Abelian integral;  bifurcation.

\vskip 0.2 true cm
\centerline{\bf { $\S$1}. Introduction and the main results}
\vskip 0.2 true cm

The determination of limit cycles is one important problem in the qualitative theory of real planar differential systems. Recently, stimulated by non-smooth phenomena in the real world such as control system [1], economics [12], nonlinear oscillations [23], and biology [13], the investigation of limit cycles for piecewise smooth differential systems has attracted many attentions.

There have been many scholars study the number of limit cycles for piecewise smooth differential systems. In [17], Llibre and Mereu studied the number of limit cycles bifurcating from the period annuluses of quadratic isochronous centers $(S_{1})$ and $(S_{2})$ when they are perturbed inside a class of piecewise smooth quadratic polynomial differential systems. J. Yang and L. Zhao improved the Theorem 1.1 and 1.2 of [17] and they obtained the sharp upper bounds of the number of limit cycles in [25]. Recently, S. Sui and L. Zhao [22] considered the bifurcation of limit cycles for generic L-V and B-T systems, and bounded the number of zeros of first order Melnikov function which controls the number of limit cycles bifurcating from the center. For more, one is recommended to see [4,18,19].

It is known that [7] proved that after an affine change of variables and a rescaling of the independent variable any cubic Hamiltonian can be transformed into the following form
$$
H(x,y)=\frac{1}{2}(x^{2}+y^{2})-\frac{1}{3}x^{3}+axy^{2}+\frac{1}{3}by^{3},
\eqno(1.1)$$
where $a$, $b$ are parameters  in the region
$$
G=\left\{(a,b):-\frac{1}{2}\leq a\leq1,0\leq b\leq (1-a)(1+2a)^{1/2}\right\}.
$$
Let
$$
X_{H}:
\left\{{\begin{aligned}
&\dot{x}=H_{y}(x,y),\\
&\dot{y}=-H_{x}(x,y),
\end{aligned}}\right.
$$
where the form of $H(x,y)$ is (1.1). Then fields $X_{H}$ are generic if $(a,b)\in G\backslash\partial G$ and degenerate if $X_{H}\in \partial G$. A lot of scholars have studied the cyclicity of period annulus of $X_{H}$ under quadratic perturbations such as [2,3,11,14,26], and it is well-known that if $(a,b)\in G\backslash(1,0)$ the cyclicity of period annulus of $X_{H}$ under quadratic perturbations equals two, and it equals three if $(a,b)\in(1,0)$. By using Abelian integra method, Ilive [8] estimated that there are at most $5n+15$ limit cycles bifurcating from the period annuluses of (1.1) with $(a,b)\in G$ under continuous perturbations of arbitrary polynomials with degree $n$.

Motivated by [8], [24] and [26], in this paper, we study the upper bound of the number of limit cycles bifurcating from the period annuluses of quadratic Hamilton system $S^{(2)}$ (two-point heteroclinic orbits) and $S^{(3)}$ (three-point heteroclinic orbits) when they are perturbed inside any discontinuous polynomial differential systems of degree $n$. There has three cases for $S^{(2)}$, elliptic segment, hyperbolic segment and parabolic segment (see Fig.\,1). $S^{(3)}$ is called Hamiltonian triangle (see Fig.\,2).

The main results are follows.

\vskip 0.1 true cm
\noindent{\bf Theorem 1.1.}  Let $0<|\varepsilon|\ll 1$, and $p^\pm(x,y)$ and $q^\pm(x,y)$ are any polynomials of degree $n$. Consider the following perturbations of the system:
$$
\left(
  \begin{array}{c}
          \dot{x}\\
          \dot{y}
   \end{array}
   \right)=\begin{cases}
   \left(
    \begin{array}{c}
        y+2axy+by^2+\varepsilon p^+(x,y)\\
        -x+x^2-ay^2+\varepsilon q^+(x,y)
   \end{array}
   \right), \quad y>0,\\
   \,
   \left(
    \begin{array}{c}
       y+2axy+by^2+\varepsilon p^-(x,y)\\
       -x+x^2-ay^2+\varepsilon q^-(x,y)
   \end{array}
   \right), \quad y<0.\\
   \end{cases}
   \eqno(1.2)$$
Then, by using the first order Melnikov function in $\varepsilon$,
 the upper bounds of the number of limit cycles of systems (1.2) bifurcating from the each period annuluses are $25n+161$($n\geq3$) and $24n+126$($n\geq3$) for
 $S^{(2)}$ and $S^{(3)}$  respectively(counting the multiplicity).

\begin{figure}[htbp]
\centering
\subfigure[elliptic segment]{
\includegraphics[width=3.8cm]{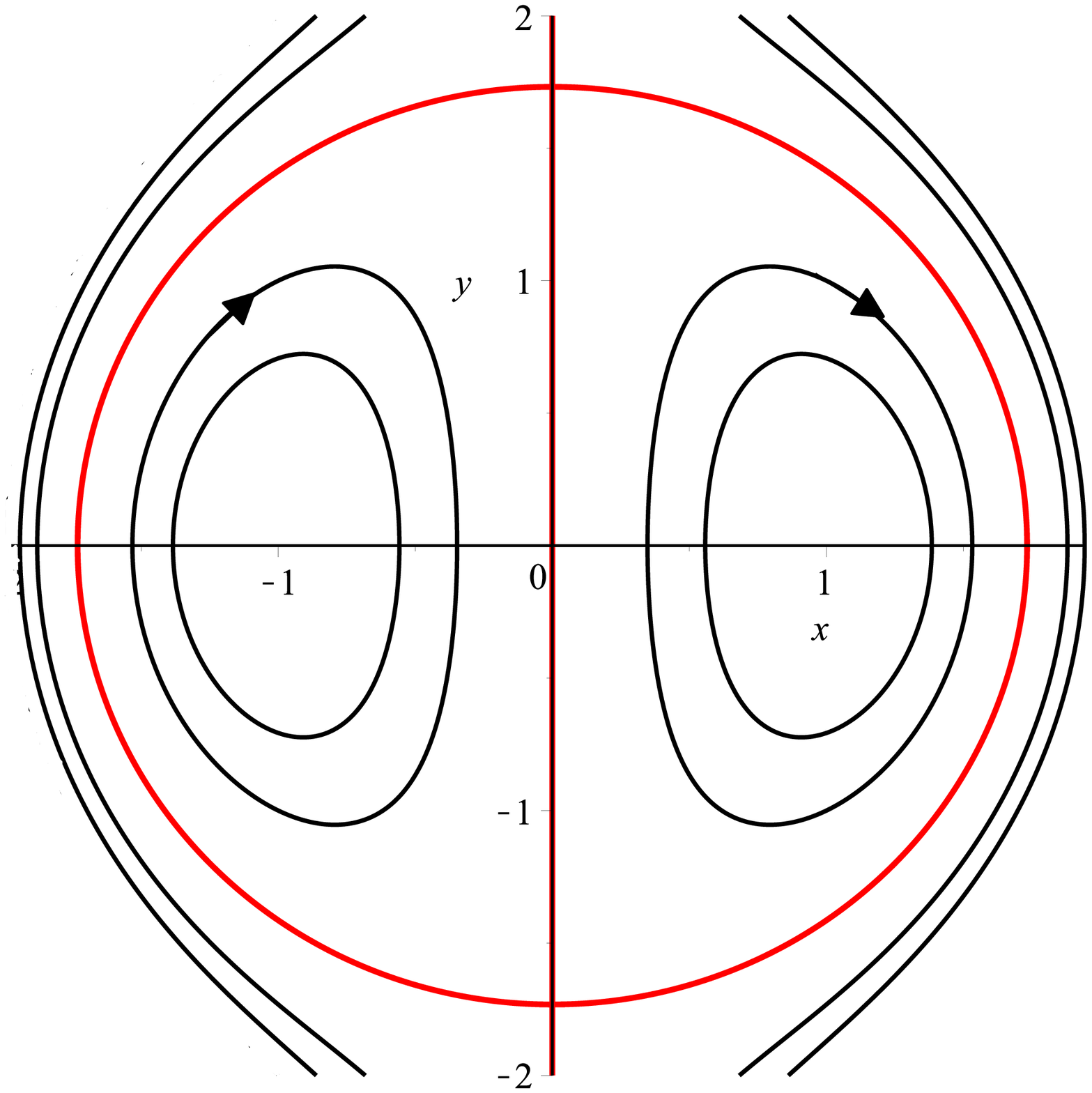}}
\quad
\subfigure[hyperbolic segment]{\includegraphics[width=3.8cm]{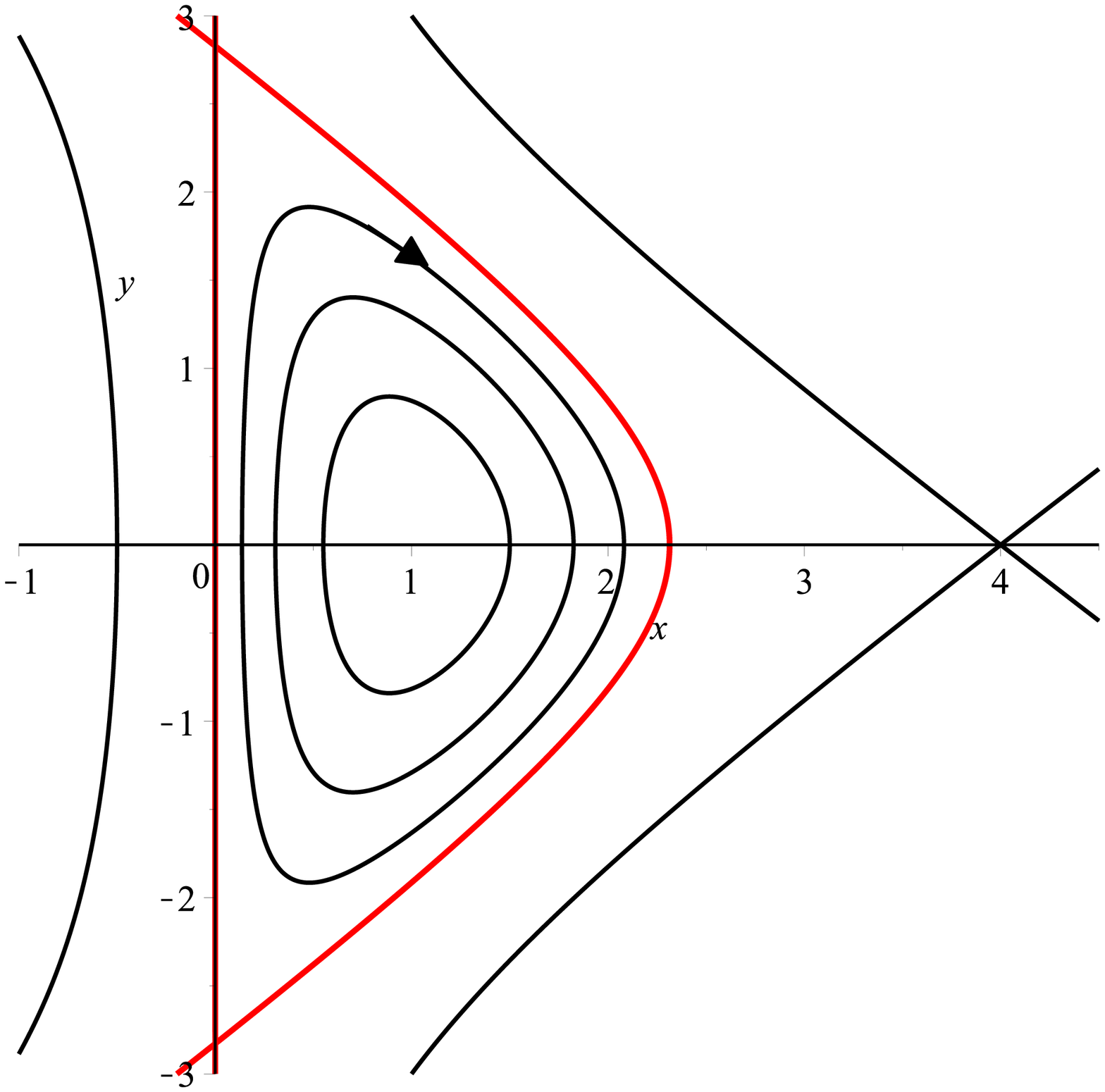}}
\quad
\subfigure[ parabolic segment]{\includegraphics[width=3.7cm]{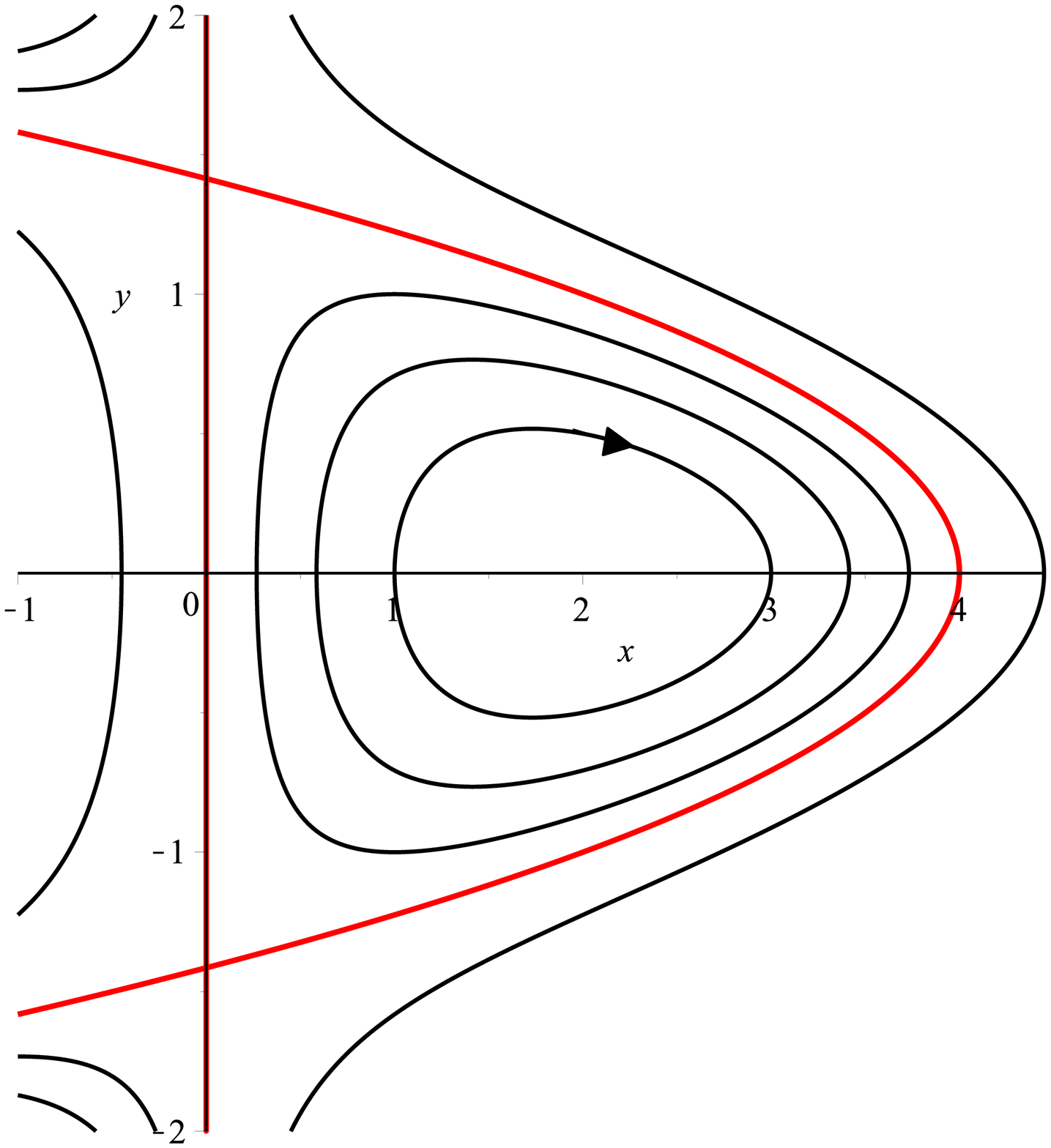}}
\caption{The phase portraits of $S^{(2)}$.
For the case of (a),  we have $a\in(-\frac{1}{2},\frac{1}{2})$ and $b=(1-a)(1+2a)^{1/2}$ in (1.2).  For the case of (b), we have $a\in(\frac{1}{2},1)$ and  $b=(1-a)(1+2a)^{1/2}$.
For the case of (c), we have $a=\frac{1}{2}$ and  $b=\frac{1}{\sqrt{2}}$. }
\end{figure}

\begin{figure}[htbp]
\centering
{\includegraphics[width=4.5cm]{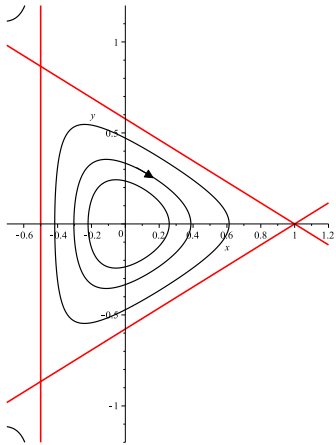}}
\caption{The phase portraits of $S^{(3)}$: $a=1$ and $b=0$.}
\end{figure}

\vskip 0.3 true cm
\noindent{\bf Remark 1.2.} (i) The techniques we use mainly include the first order Melnikov function, Picard-Fuchs equation and Chebyshev criterion. By [9], we know that the number of zeros of the first order Melnikov function $M(h)$ controls the number of limit cycles of systems (1.2) if $M(h)\neq0$ in the corresponding period annulus.

\noindent(ii) For the parabolic segment of $S^{(2)}$, the result is meaningful for $n\geq2$, $S^{(3)}$ and the rest cases of $S^{(2)}$ is meaningful for $n\geq3$.
\vskip 0.1 true cm

The paper is organized as follows. In $\S$2, we will give some preliminaries. In $\S3$ and $\S4$, we will prove Theorem 1.1. First, we will obtain the algebraic structure of the first order Melnikov functions $M(h)$ for $h\in\Sigma$, which are more complicated than the Melnikov function corresponding to the continuous perturbations. Then we prove that there exists a second-order differential operator that can simplify $M(h)$. Finally, the main results are proved by using the Chebyshev space. Noting that there has six generators for elliptic segment, hyperbolic segment and Hamilton triangle, so that we use some different techniques to reduce the number of generators.

Throughout the paper, we denote by $\#\{f(h)=0,h\in(s,t)\}$ the number of isolated zeros of $f(h)$ on$(s,t)$ taking into account the multiplicity, and denote by $A^{T}$ the transpose of matrix $A$. We will give clear instructions if $P_{k}(h)$ express the polynomials of degree at most $k$, the others do not have that meaning.

\vskip 0.2 true cm
\centerline{\bf { $\S$2.} Preliminaries}
\vskip 0.2 true cm

We first introduce the first order Melnikov function of discontinuous differential systems. Consider the following systems:
$$
(\dot{x},\ \dot{y})=\begin{cases}
       (P^+(x,y)+\varepsilon p^+(x,y),Q^+(x,y)+\varepsilon q^+(x,y)),\ \ y>0,\\
       (P^-(x,y)+\varepsilon p^-(x,y),Q^-(x,y)+\varepsilon q^-(x,y)),\ \ y<0,
\end{cases}
\eqno(2.1)$$
where $0<|\varepsilon|\ll 1$, and $p^\pm(x,y)$ and $q^\pm(x,y)$ are polynomials with degree $n$. System $(2.1)$ has two subsystems:
$$
\left\{{\begin{aligned}
\dot{x}=P^{+}(x,y)+\varepsilon p^{+}(x,y),\\
\dot{y}=Q^{+}(x,y)+\varepsilon q^{+}(x,y),
\end{aligned}}~~~~~~y>0,
\right.\eqno(2.2)$$
and
$$
\left\{{\begin{aligned}
\dot{x}=P^{+}(x,y)+\varepsilon p^{-}(x,y),\\
\dot{y}=Q^{+}(x,y)+\varepsilon q^{-}(x,y),
\end{aligned}}~~~~~~y<0.
\right.\eqno(2.3)$$
Suppose that system $(2.2)_{\varepsilon=0}$ is integrable with the first integral $H^{+}(x,y)$ and integrating factor $\mu_{1}(x,y)$, and system $(2.3)_{\varepsilon=0}$ is integrable with the first integral $H^{-}(x,y)$ and integrating factor $\mu_{2}(x,y)$. We also suppose that $(2.1)_{\varepsilon=0}$ has a family of periodic orbits around the origin and satisfies the following two assumptions.

\vskip 0.3 true cm

{\bf Assumption (I).} There exist an interval $\Sigma=(\alpha, \beta)$, and two points $A(h)=(a(h),0)$ and $B(h)=(b(h),0)$ such that for $h\in{\Sigma}$
$$
H^{+}(A(h))=H^{+}(B(h))=h,~~H^{-}(A(h))=H^{-}(B(h))=\tilde{h},~~a(h)\neq b(h).$$

{\bf Assumption (II).} The subsystem $(2.2)_{\varepsilon=0}$ has an orbital arc $L_{h}^{+}$ starting from $A(h)$ and ending at $B(h)$ defined by $H^{+}(x,y)=h$, $y\geq0$. The subsystem $(2.3)_{\varepsilon=0}$ has an orbital arc $L_{h}^{-}$ starting from $B(h)$ and ending at $A(h)$ defined by $H^{-}(x,y)=\tilde{h}$, $y<0$.

\vskip 0.3 true cm

Under the Assumptions (I) and (II), $(2.1)_{\varepsilon=0}$ has a family of non-smooth periodic orbits $L_{h}=L_{h}^{+}\cup L_{h}^{-}(h\in \Sigma)$. For definiteness, we assume that the orbits $L_{h}$ for $h\in \Sigma$ orientate clockwise(see Fig.\,3). The authors [10] established a bifurcation function $F(h,\varepsilon)$ for $(2.1)$. Let $F(h,0)=M(h)$. In [9] and [21], the authors obtained the following results.

\begin{figure*}[!!ht]
\centering
{\includegraphics[scale=0.3]{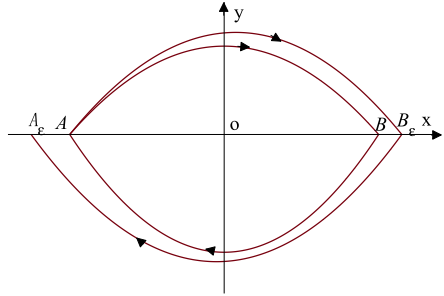}}
\caption{The Poincar\'{e} map related to $y=0$.}
\end{figure*}

\vskip 0.3 true cm
\noindent{\bf Lemma 2.1.}([9,21]). Under the assumptions (I) and (II), we have

(i) If $M(h)$ has $k$ zeros in $h$ on the interval $\Sigma$ with each having an odd multiplicity, then $(2.2)$ has at least $k$ limit cycles bifurcating from the period annulus for $0<\left|\varepsilon\right|\ll1$.

(ii) If $M(h)$ has at most $k$ zeros in $h$ on the interval $\Sigma$, taking into account the multiplicity, then there exist at most $k$ limit cycles of $(2.1)$ bifurcating from the period annulus.

(iii) The first order Melnikov function $M(h)$ of system $(2.1)$ has the following form
$$
\begin{aligned}
M(h)&=\frac{H_{x}^{+}(A)}{H_{x}^{-}(A)}\left[\frac{H_{x}^{-}(B)}{H_{x}^{+}(B)}\int_{L_{h}^{+}}\mu_{1}(x,y)q^{+}(x,y)dx-\mu_{1}(x,y)p^{+}(x,y)dy\right.\notag\\
&\left.+\int_{L_{h}^{-}}\mu_{2}(x,y)q^{-}(x,y)dx-\mu_{2}(x,y)p^{-}(x,y)dy\right].
\end{aligned}
\eqno(2.4)$$

Let $\mathbb{V}$ be a finite-dimensional vector space of functions, real-analytic on an open interval $\mathbb{I}$. Next, we give the relation of the number of zeros about second order linear homogeneous equation and non-homogeneous equation which will be used in the proof of Theorem 1.1.

\vskip 0.3 true cm
\noindent{\bf Definition 2.1.}([5]). We say that $\mathbb{V}$ is a {\bf Chebyshev space}, provided that each non-zero function in $\mathbb{V}$ has at most $dim(\mathbb{V})-1$ zeros, counted with multiplicity.

Let $\mathbb{S}$ be the solution space of a second order linear analytic differential equation
$$
x^{''}+a_{1}(t)x^{'}+a_{2}(t)x=0
\eqno(2.5)$$
on an open interval $\mathbb{I}$.

\vskip 0.1 true cm
\noindent{\bf Lemma 2.3.}([5]). The solution space $\mathbb{S}$ of $(2.5)$ is a Chebyshev space of the interval $\mathbb{I}$ if and only if there exists a nowhere vanishing solution $x_{0}(t)\in \mathbb{S}$($x_{0}(t)\neq0, \forall t\in \mathbb{I}$).

\vskip 0.1 true cm
\noindent{\bf Lemma 2.4.}([5]). Suppose the solution space of the homogeneous equation $(2.5)$ is a Chebyshev space and let $R(t)$ be an analytic function on $\mathbb{I}$ having $l$ zeros (counted with multiplicity). Then every solution $x(t)$ of the non-homogeneous equation
$$x^{''}+a_{1}(t)x^{'}+a_{2}(t)x=R(t)$$
has at most $l+2$ zeros on $\mathbb{I}$.

\vskip 0.1 true cm
\noindent{\bf Lemma 2.5.}([6]). Consider the following function:
$$
F(x)=P^{0}(x)+\sum_{j=1}^{k}P^{j}(x)\frac{1}{\sqrt{x+c_{j}}},
$$
where $P^{j}(x)(j=0,1,...,k)$ is a real polynomial function, $c_{j}(j=0,1,...,k)$ is real constants. $Z(F)$ is the number of zeros of $F(x)$ on $({\rm max}_{j=1,2,...,k}\left\{-c_{j}\right\}, +\infty)$ (taking into account the multiplicity), then
$$
Z(F)\leq k({\rm max}_{j=1,2,...,k}\left\{{\rm deg}(P^{j})\right\}+1)+{\rm deg}(P^{0}),~~{\rm deg}(0)=-1.
$$

\vskip 0.3 true cm
\centerline{\bf { $\S$3.} Proof of the results on  $S^{(2)}$}
\vskip 0.2 true cm

As the $a$, $b$ varying in value, the singularity of $S^{(2)}$ has changed, therefore we can classify $S^{(2)}$ as follows:

(1) elliptic segment: $a\in(-\frac{1}{2},\frac{1}{2})$, $b=(1-a)(1+2a)^{1/2}$;

(2) hyperbolic segment: $a\in(\frac{1}{2},1)$, $b=(1-a)(1+2a)^{1/2}$;

(3) parabolic segment: $a=\frac{1}{2}$, $b=(1-a)(1+2a)^{1/2}$.

In this section, we will mainly prove the elliptic segment, the others are similar. And the same symbol means different significance in each case in order to make reading and writing easier.

\vskip 0.3 true cm
{ {\bf { $\S$3.1.} The case of the elliptic segment.}}
\vskip 0.3 true cm
In this section, let
$$
x_{1}=\frac{b}{1-a}y-x+1,~~y_{1}=\frac{1-a}{b}\sqrt{3(2-\lambda)}y+\sqrt{3(2-\lambda)}x,~~\lambda=\frac{1-a(1+2a)}{1+a}.
$$
Denote $\Sigma_{1}=(\lambda-3,0)$ and $\Sigma_{2}=(0,(\lambda-2)^{2}(\lambda+1)/\lambda^{2})$, yhen (1.1) and (1.2) transform into
$$
H(x,y)=xy^{2}+\lambda x^{3}-3(\lambda-1)x^{2}+3(\lambda-2)x=h,~h\in\Sigma_{1}\cup\Sigma_{2}
\eqno(3.1)$$
$$
\left(
  \begin{array}{c}
          \dot{x}\\
          \dot{y}
   \end{array}
   \right)=\begin{cases}
   \left(
    \begin{array}{c}
        2xy+\varepsilon p^+(x,y)\\
        -y^{2}-3\lambda x^{2}+6(\lambda-1)x-3(\lambda-2)+\varepsilon q^+(x,y)
   \end{array}
   \right), \quad y>0,\\
   \left(
    \begin{array}{c}
       2xy+\varepsilon p^-(x,y)\\
       -y^{2}-3\lambda x^{2}+6(\lambda-1)x-3(\lambda-2)+\varepsilon q^-(x,y)
   \end{array}
   \right), \quad y<0.\\
   \end{cases}
   \eqno(3.2)$$
where
$$p_{n}^{\pm}(x,y)=\sum\limits_{i+j=0}^na_{i,j}^{\pm}x^{i}y^{j},~~q_{n}^{\pm}(x,y)=\sum\limits_{i+j=0}^nb_{i,j}^{\pm}x^{i}y^{j}$$
are any polynomials of degree $n$. Here and below we shall omit the subscript 1.

For $\lambda\in(0,2)$, System (3.2) has two elementary centers $(1,0)$ and $(\lambda-2)/\lambda,0)$ corresponding to $h=\lambda-3$ and $h=(\lambda-2)^{2}(\lambda+1)/\lambda^{2}$ respectively, two saddles $(0,\pm\sqrt{3(2-\lambda)})$ corresponding to $h=0$ (see Fig.\,2(a)).

Then there exist period annulus in right half plane if $h\in\Sigma_{1}$ and in left half plane if $h\in\Sigma_{2}$, we will prove the conclusion for $h\in\Sigma_{1}$ in the following and it is similar for $h\in\Sigma_{2}$.

By Lemma 2.1, we have
$$M(h)=\int_{L_{h}^{+}}q_{n}^{+}(x,y)dx-p_{n}^{+}(x,y)dy+\int_{L_{h}^{-}}q_{n}^{-}(x,y)dx-p_{n}^{-}(x,y)dy,~h\in\Sigma_{1},\eqno(3.3)$$
where
$$
L_{h}^{+}(L_{h}^{-})=\{(x,y)|H(x,y)=h,y>0(y<0),h\in\Sigma_{1}\}.$$
For $i,~j>0$, denote
$$
J_{i,j}(h)=\int_{L_{h}^{+}}x^iy^jdx.
$$
\vskip 0.3 true cm
\noindent{\bf Lemma 3.1.} Consider system (3.2), if $n\geq3$, for $h\in\Sigma_{1}$, $M(h)$ can be expressed as
$$
M(h)={\bf \sigma_{1}(h)}{\bf U_{1}(h)}+{\bf \sigma_{2}(h)}{\bf U_{2}(h)},
\eqno(3.4)$$
where
$$
{\bf \sigma_{1}(h)}=(\alpha(h),\beta(h),\gamma(h)),~~{\bf \sigma_{2}(h)}=(\eta(h),\xi(h),\zeta(h)),
$$
$$
{\bf U_{1}(h)}=(J_{0,0}(h),J_{1,0}(h),J_{0,2}(h))^{T},~~{\bf U_{2}(h)}=(J_{0,1}(h),J_{1,1}(h),J_{2,1}(h))^{T},
$$
and $\alpha(h)$, $\beta(h)$, $\gamma(h)$, $\eta(h)$, $\xi(h)$, $\zeta(h)$ are polynomials of $h$ with
$$
{\rm deg}\alpha(h)\leq\left[\frac{n}{3}\right],~~{\rm deg}\beta(h)\leq\left[\frac{n-1}{3}\right],~~{\rm deg}\gamma(h)\leq\left[\frac{n-2}{3}\right],
$$
$$
{\rm deg}\eta(h)\leq\left[\frac{n-1}{3}\right],~~{\rm deg}\xi(h)\leq\left[\frac{n-2}{3}\right],~~{\rm deg}\zeta(h)\leq\left[\frac{n-3}{3}\right].
$$

\vskip 0.1 true cm
\noindent{\bf Proof.} For system (3.2), we assume that
$$
M(h)=\sum\limits_{\begin{subarray}{c}i+j=0\\i\geq0,j\geq0\end{subarray}}^{n}\rho_{i,j}J_{i,j}(h),\eqno(3.5)
$$
where $\rho_{i,j}$ are arbitrary real constants. In fact, Suppose that orbit $L_{h}^{+}(L_{h}^{-})$ intersects the $x$-axis at points $A_{h}(x_{A}(h),0)$ and $B_{h}(x_{B}(h),0)$ for $h\in\Sigma_{1}$. Let $D^{+}$ be the interior of $L_{h}^{+}\cup\overrightarrow{B_{h}A_{h}}$. Then by direct computation,
$$\int_{L_{h}^{+}}x^iy^jdy=\int_{L_{h}^{+}\cup\overrightarrow{B_{h}A_{h}}}x^iy^jdy-\int_{\overrightarrow{B_{h}A_{h}}\,}x^iy^jdy
=-i\iint\limits_{D}x^{i-1}y^jdxdy,$$
and
$$\int_{L_{h}^{+}}x^{i-1}y^{j+1}dx=\int_{L_{h}^{+}\cup\overrightarrow{B_{h}A_{h}}}x^{i-1}y^{j+1}dx-\int_{\overrightarrow{B_{h}A_{h}}}x^{i-1}y^{j+1}dx=
(j+1)\iint\limits_{D}x^{i-1}y^{j}dxdy.$$
Hance we have
$$\int_{L_{h}^{+}}x^iy^jdy=-\frac{i}{j+1}\int_{L_{h}^{+}}x^{i-1}y^{j+1}dx.$$
Similarly, we have
$$\int_{L_{h}^{-}}x^iy^jdy=-\frac{i}{j+1}\int_{L_{h}^{-}}x^{i-1}y^{j+1}dx.$$
Noting the symmetry of $H(x,y)$, we have
$$\int_{L_{h}^{-}}x^iy^jdx=(-1)^{j+1}\int_{L_{h}^{+}}x^iy^jdx.$$
therefore, we have
$$
\begin{aligned}
M(h)
&=\int_{L_{h}^{+}}\sum\limits_{i+j=0}^nb_{i,j}^{+}x^iy^jdx+\sum\limits_{i+j=0}^na_{i,j}^{+}\frac{i}{j+1}x^{i-1}y^{j+1}dx\\
&+\int_{L_{h}^{-}}\sum\limits_{i+j=0}^nb_{i,j}^{-}x^iy^jdx+\sum\limits_{i+j=0}^na_{i,j}^{-}\frac{i}{j+1}x^{i-1}y^{j+1}dx\\
&=\int_{L_{h}^{+}}\sum\limits_{i+j=0}^nb_{i,j}^{+}x^iy^jdx+\sum\limits_{i+j=0}^na_{i,j}^{+}\frac{i}{j+1}x^{i-1}y^{j+1}dx\\
&+(-1)^{j+1}\int_{L_{h}^{+}}\sum\limits_{i+j=0}^nb_{i,j}^{-}x^iy^jdx+\sum\limits_{i+j=0}^na_{i,j}^{-}\frac{i}{j+1}x^{i-1}y^{j+1}dx,
\end{aligned}
$$
than we can get (3.5).
Differentiating (3.1) with respect to $x$, we obtain
$$
y^{2}+3\lambda x^{2}-6(\lambda-1)x+3(\lambda-2)=0.
\eqno(3.6)$$
Multiplying (3.1) and (3.6) by $x^{i}y^{j}dx$, integrating over $L_{h}^{+}$, we have
$$
\lambda J_{i,j}(h)+J_{i-2,j+2}=hJ_{i-3,j}(h)+3(\lambda-1)J_{i-1,j}(h)-3(\lambda-2)J_{i-2,j}(h),
\eqno(3.7)$$
$$
\frac{j-2i-2}{j}J_{i,j}(h)+3\lambda J_{i+2,j-2}(h)=6(\lambda-1)J_{i+1,j-2}(h)-3(\lambda-2)J_{i,j-2}(h).
\eqno(3.8)$$
Elementary manipulations reduce Eqs. (3.7) and (3.8) to
$$
\frac{2i+2j+2}{j}J_{i,j}(h)=3hJ_{i-1,j-2}(h)+3(\lambda-1)J_{i+1,j-2}(h)-6(\lambda-2)J_{i,j-2}(h),
\eqno(3.9)$$
$$
\begin{aligned}
\frac{2i+2j+2}{j+2}\lambda J_{i,j}(h)&=\frac{2i-j-4}{j+2}hJ_{i-3,j}(h)\\
&+3(\lambda-1)\frac{2i+j}{j+2}J_{i-1,j}(h)-3(\lambda-2)\frac{2i-2}{j+2}J_{i-2,j}(h).
\end{aligned}
\eqno(3.10)$$
Let $i=0$ in (3.8), we have
$$
\frac{j-2}{j}J_{0,j}(h)=6(\lambda-1)J_{0,j-2}(h)-3(\lambda-2)J_{0,j-2}(h)-3\lambda J_{2,j-2}(h).
\eqno(3.11)$$
We will prove the conclusion by induction on $n$. When $n=3$, (3.9)-(3.11) give
$$
\begin{aligned}
J_{3.0}(h)&=(\frac{1}{4\lambda}h-\frac{9}{4}+\frac{27}{4\lambda}-\frac{9}{2\lambda^{2}})J_{0,0}(h)+(3-\frac{6}{\lambda}+\frac{9}{2\lambda^{2}})J_{1,0}(h),\\
J_{1,2}(h)&=(\frac{3}{4}h-\frac{3}{4}\lambda+\frac{9}{4}-\frac{3}{2\lambda})J_{0,0}(h)+\frac{3}{2\lambda}J_{1,0}(h),\\
J_{0,3}(h)&=18(\lambda-1)J_{1,1}(h)-9(\lambda-2)J_{0,1}(h)-9\lambda J_{2,1}(h),
\end{aligned}
$$
which yields the result for $n=3$. Suppose that the result holds for $i+j\leq k(k\geq3)$. Then for $i+j=k+1(k\geq2)$, taking $(i,j)=(k+1,0),(k,1)$ in (3.10), $(i,j)=(k-1,2),...,(2,k-1),(1,k)$ in (3.9) and $(i,j)=(0,k+1)$ in (3.11) respectively, we can obtain that
$$
{\bf G}\left(\begin{matrix}
          J_{0,k+1}(h)\\
          J_{1,k}(h)\\
          J_{2,k-1}(h)\\
          \vdots\\
          J_{k-1,2}(h)\\
          J_{k,1}(h)\\
          J_{k+1,0}(h)
          \end{matrix}\right)\ \
=\left(\begin{matrix}
        C(h)\\
        3hJ_{0,k-2}(h)+3(\lambda-1)J_{2,k-2}(h)-6(\lambda-2)J_{1,k-2}(h)\\
        3hJ_{1,k-3}(h)+3(\lambda-1)J_{3,k-3}(h)-6(\lambda-2)J_{2,k-3}(h)\\
        \vdots\\
        3hJ_{k-2,0}(h)+3(\lambda-1)J_{k,0}(h)-6(\lambda-2)J_{k-1,0}(h)\\
        h\frac{2k-5}{3}J_{k-3,1}(h)-3(\lambda-1)(2k-5)J_{k-1,1}(h)-3(\lambda-2)(k-1)J_{k-2,1}(h)\\
        h(k-1)J_{k-2,0}(h)-3(\lambda-1)(k+1)J_{k,0}(h)-3(\lambda-2)kJ_{k-1,0}(h)
        \end{matrix}\right),
$$
where
$$
\begin{aligned}
C(h)=&6(\lambda-1)J_{1,k-1}(h)-3(\lambda-2)J_{0,k-1}(h)-9h\lambda\frac{k-1}{4+2k}J_{1,k-3}(h)\\
&-9\lambda(\lambda-1)\frac{k-1}{4+2k}J_{3,k-3}(h)+18\lambda(\lambda-2)\frac{k-1}{4+2k}J_{2,k-3}(h),
\end{aligned}
$$
$$
{\bf G}=\left(\begin{matrix}
       &\frac{k-1}{k+1}~~&0~~&0~~&\dots~~&0~~&0~~&0\\
       &0~~&\frac{2k+4}{k}~~&0~~&\dots~~&0~~&0~~&0\\
       &0~~&0~~&\frac{2k+4}{k-1}~~&\dots~~&0~~&0~~&0\\
       &\vdots~~&\vdots~~&\vdots~~&\vdots~~&\vdots~~&\vdots\\
       &0~~&0~~&0~~&\dots~~&k+2~~&0~~&0\\
       &0~~&0~~&0~~&\dots~~&0~~&\frac{2k+4}{3}~~&0\\
       &0~~&0~~&0~~&\dots~~&0~~&0~~&k+2
       \end{matrix}\right),
$$
and ${\rm det}{\bf G}\neq0$.
By the induction hypothesis we obtain the expression (3.4). Next we estimate the degree of polynomials of $\alpha(h)-\zeta(h)$, taking $J_{1,k}(h)$ as an example.
$$
\begin{aligned}
J_{1,k}(h)=&h(\alpha_{\left[\frac{k-2}{3}\right]}J_{0,0}+\beta_{\left[\frac{k-3}{3}\right]}J_{1,0}+\gamma_{\left[\frac{k-4}{3}\right]}J_{0,2}
+\eta_{\left[\frac{k-3}{3}\right]}J_{0,1}+\xi_{\left[\frac{k-4}{3}\right]}J_{1,1}+\zeta_{\left[\frac{k-5}{3}\right]}J_{2,1})\\
&+(\alpha_{\left[\frac{k}{3}\right]}J_{0,0}+\beta_{\left[\frac{k-1}{3}\right]}J_{1,0}+\gamma_{\left[\frac{k-2}{3}\right]}J_{0,2}
+\eta_{\left[\frac{k-1}{3}\right]}J_{0,1}+\xi_{\left[\frac{k-2}{3}\right]}J_{1,1}+\zeta_{\left[\frac{k-3}{3}\right]}J_{2,1})\\
&+\alpha_{\left[\frac{k-1}{3}\right]}J_{0,0}+\beta_{\left[\frac{k-2}{3}\right]}J_{1,0}+\gamma_{\left[\frac{k-3}{3}\right]}J_{0,2}
+\eta_{\left[\frac{k-2}{3}\right]}J_{0,1}+\xi_{\left[\frac{k-3}{3}\right]}J_{1,1}+\zeta_{\left[\frac{k-4}{3}\right]}J_{2,1}\\
&:=\tilde{\alpha}J_{0,0}+\tilde{\beta}J_{1,0}+\tilde{\gamma}J_{0,2}+\tilde{\eta}J_{0,1}+\tilde{\xi}J_{1,1}+\tilde{\zeta}J_{2,1},
\end{aligned}
$$
where $\alpha_{m}-\zeta_{m}$ represent the polynomials of h satisfying ${\rm deg}\alpha_{m}\leq m$, and $\beta_{m}-\zeta_{m}$ are the same.
It is easy to check
$$
{\rm deg}\tilde{\alpha}\leq {\rm max}\left\{\left[\frac{k-1}{3}\right]+1,~\left[\frac{k}{3}\right]\right\}=\left[\frac{k+1}{3}\right].
$$
In the similar way, we can end the proof. $\diamondsuit$

\vskip 0.3 true cm
\noindent{\bf Lemma 3.2.}
If $n\geq3$, for $h\in\Sigma_{1}$,

\noindent(1) the vector functions ${\bf U_{1}(h)}$ and ${\bf U_{2}(h)}$ satisfy respectively the following Picard-Fuchs equations:
$$
{\bf U_{1}(h)}=(B_{1}h+C_{1}){\bf U_{1}^{'}(h)},
\eqno(3.12)$$
$$
{\bf U_{2}(h)}=(B_{2}h+C_{2}){\bf U_{2}^{'}(h)},
\eqno(3.13)$$
where
$$
B_{1}=\left(\begin{matrix}
             &3~~&0~~&0\\
             &\frac{3}{2}-\frac{3}{2\lambda}~~&\frac{3}{2}~~&0\\
             &-\frac{33}{8}\lambda+\frac{33}{4}+\frac{15}{8\lambda}~~&\frac{15}{8}\lambda-\frac{15}{8}~~&1
      \end{matrix}\right),
$$
$$
C_{1}=\left(\begin{matrix}
            &9-3\lambda-\frac{6}{\lambda}~~&\frac{6}{\lambda}~~&0\\
            &-\frac{3}{2}\lambda+6-\frac{21}{2\lambda}+\frac{9}{\lambda^{2}}~~&-\frac{3}{2}\lambda+\frac{9}{2}+\frac{6}{\lambda}-\frac{9}{\lambda^{2}}~~&0\\
            &\frac{33}{8}\lambda^{2}-\frac{165}{8}\lambda+\frac{219}{8}+\frac{3}{8\lambda}-\frac{90}{8\lambda^{2}}~~&-\frac{15}{8}\lambda^{2}+\frac{30}{4}\lambda-\frac{81}{8}+\frac{21}{8\lambda}+\frac{90}{8\lambda^{2}}~~&0
      \end{matrix}\right),
$$
$$
B_{2}=\left(\begin{matrix}
            &\frac{3}{2}~~&0~~&0\\
            &\frac{8}{5}(1-\frac{1}{\lambda})~~&1~~&0\\
            &\frac{3}{320\lambda^{2}}(-179\lambda^{2}+358\lambda+77)~~&\frac{3}{8}(1-\frac{1}{\lambda})~~&\frac{3}{4}
      \end{matrix}\right),
$$
$$
C_{2}=\left(\begin{matrix}
            &0~~&-3\lambda+6~~&\frac{3}{2}(\lambda-1)\\
            &0~~&-\frac{21}{5}\lambda+\frac{63}{5}-\frac{42}{5\lambda}~~&\frac{8}{5}\lambda-\frac{16}{5}+\frac{18}{5\lambda}\\
            &0~~&\frac{3}{320\lambda^{2}}(318\lambda^{3}-1272\lambda^{2}+918\lambda+708)~~&\frac{3}{320\lambda^{2}}(-259\lambda^{3}+777\lambda^{2}-41\lambda-477)
      \end{matrix}\right).
$$
And we have
$$
D_{1}(h):={\rm det}(B_{1}h+C_{1})=\frac{9}{2\lambda^{2}}h(\lambda-h-3)(\lambda^{3}-\lambda^{2}h-3\lambda^{2}+4),
\eqno(3.14)$$
$$
D_{2}(h):={\rm det}(B_{1}h+C_{1})=\frac{9}{8\lambda^{2}}h(\lambda-h-3)(\lambda^{3}-\lambda^{2}h-3\lambda^{2}+4).
\eqno(3.15)$$
\noindent{(2)}
$$
{\rm det}(B_{1}h+C_{1}){\bf U_{1}^{''}(h)}=
        \left(\begin{matrix}
          &g_{11}(h)~~&g_{12}(h)\\
          &g_{21}(h)~~&g_{22}(h)\\
          &g_{31}(h)~~&g_{32}(h)
        \end{matrix}\right)
        \left(\begin{matrix}
        J_{0,0}^{'}(h)\\
        J_{1,0}^{'}(h)
        \end{matrix}\right),
\eqno(3.16)$$
$$
{\rm det}(B_{2}h+C_{2})\left(\begin{matrix}
                         J_{0,1}^{''}(h)\\
                         J_{1,1}^{''}(h)\\
                         Z^{''}(h)
                   \end{matrix}\right)
        =\left(\begin{matrix}
        &k_{11}(h)~~&k_{12}(h)\\
        &k_{21}(h)~~&k_{22}(h)\\
        &k_{31}(h)~~&k_{32}(h)
        \end{matrix}\right)
        \left(\begin{matrix}
        J_{0,1}^{'}(h)\\
        Z^{'}(h)
        \end{matrix}\right),
\eqno(3.17)$$
where
$$
Z(h)=\frac{3}{8}(\frac{1}{\lambda}-1)J_{1,1}(h)+\frac{1}{4}J_{2,1}(h),
\eqno(3.18)$$
and
$$
\begin{aligned}
g_{11}(h)&=3h[\lambda^{3}-\lambda^{2}(h+3)-\lambda+3]/\lambda^{2},~~~~~~~~~~~~~~g_{12}(h)=3h/\lambda,\\
g_{21}(h)&=3h[\lambda^{3}-(h+4)\lambda^{2}+(h+1)\lambda+6]/ 2\lambda^{2},~~~g_{22}(h)=3[\lambda^{2}-(h+3)\lambda+2]/ 2\lambda,\\
g_{31}(h)&=9[\lambda^{5}-(2h+8)\lambda^{4}+(h^{2}+21)\lambda^{3}+(2h^{2}+12h+14)\lambda^{2}-(4h+20)\lambda+8h+24]/\lambda^{2},\\
g_{32}(h)&=-45[\lambda^{5}-(7+2h)\lambda^{4}+(h^{2}-6h+15)\lambda^{3}-(h^{2}+5)\lambda^{2}-(4h+16)\lambda+4h+12]/ 8\lambda^{2},\\
k_{11}(h)&=3[459\lambda^{5}-(419h+2754)\lambda^{4}-(40h^{2}-1257h-4374)\lambda^{3}-(985h-864)\lambda^{2}\\
&+(147h-5265)\lambda+1458]/320\lambda^{2},\\
k_{12}(h)&=3[\lambda^{3}-(4+h)\lambda^{2}+(1+h)\lambda+6]/2\lambda,\\
k_{21}(h)&=3h[64\lambda^{4}+(-64h-256)\lambda^{3}+(64h-179)\lambda^{2}+870\lambda-243]/160\lambda^{3},~~~k_{22}(h)=3h/\lambda,\\
k_{31}(h)&=3h[435\lambda^{5}-(+435h+2175)\lambda^{4}+(870h+2418)\lambda^{3}-(179h-1446)\lambda^{2}-2853\lambda\\
&+729]/1280\lambda^{4},\\
k_{32}(h)&=-3h(\lambda^{3}-\lambda^{2}h-3\lambda)/8\lambda^{2}.
\end{aligned}
$$

\vskip 0.1 true cm
\noindent{\bf Proof.}
By direct computation, for $h\in\Sigma_{1}$, we have
$$
J_{i,j}^{'}(h)=\int_{x_{A}(h)}^{x_{B}(h)}jx^{i}y^{j-1}\frac{\partial{y}}{\partial{h}}dx+x_{B}(h)^{i}y(x_{B}(h),h)^{j}\frac{\partial{x_{B}(h)}}{\partial{h}}
-x_{A}(h)^{i}y(x_{A}(h),h)^{j}\frac{\partial{x_{A}(h)}}{\partial{h}},
$$
$L_{h}^{+}$ intersects the right $x$-axis at points $(x_{A}(h),0)$ and $(x_{B}(h),0)$.
Differentiating $H(x,0)=h$ with respect to $x$, we have
$$
\frac{\partial{x}}{\partial{h}}=\frac{1}{(x-1)(x-(\lambda-2)/\lambda)},
$$
Through the analysis of the singularity of the system, we can obtain $x_{A}(h)\neq1$, $x_{A}(h)\neq(\lambda-2)/\lambda$ and $x_{B}(h)\neq1$, $x_{B}(h)\neq(\lambda-2)/\lambda$. Then we have
$$
\frac{\partial{x_{A}(h)}}{\partial{h}}\neq\infty,~~~\frac{\partial{x_{B}(h)}}{\partial{h}}\neq\infty.
$$
Hance by $y(x_{A}(h),h)^{j}=y(x_{B}(h),h)^{j}=0$, we can obtain that
$$
J_{i,j}^{'}(h)=j\int_{x_{A}(h)}^{x_{B}(h)}x^{i}y^{j-1}\frac{\partial{y}}{\partial{h}}dx.
$$
Differentiating (3.1) with respect to $h$, we have
$$
\frac{\partial{y}}{\partial{h}}=\frac{1}{2xy},
$$
then
$$
J_{i,j}(h)=\frac{2}{j+2}J_{i+1,j+2}^{'}(h).
$$
By (3.9), we have
$$
\frac{2i+2j+8}{j+2}J_{i+1,j+2}(h)=3hJ_{i,j}(h)+3(\lambda-1)J_{i+2,j}(h)-6(\lambda-2)J_{i+1,j}(h),
$$
then we can obtain
$$
(i+j+1)J_{i,j}(h)=3hJ_{i,j}^{'}(h)+3(\lambda-1)J_{i+2,j}^{'}(h)-6(\lambda-2)J_{i+1,j}^{'}(h).
\eqno(3.19)$$
Combining (3.9), (3.10), (3.11) and (3.19), we can obtain (3.12) and (3.13).

Differentiating both side (3.12) and (3.13) yields,
$$
{\rm det}(B_{1}h+C_{1}){\bf U_{1}^{''}(h)}=(B_{1}h+C_{1})^{*}(E-B_{1}){\bf U_{1}^{'}(h)},
$$
$$
{\rm det}(B_{2}h+C_{2}){\bf U_{2}^{''}(h)}=(B_{2}h+C_{2})^{*}(E-B_{2}){\bf U_{2}^{'}(h)},
\eqno(3.20)$$
where $E$ is $3\times3$ identity matrix.
substituting (3.18) into (3.20) gives (3.16) and (3.17).
This ends the proof. $\diamondsuit$

\vskip 0.3 true cm
\noindent{\bf Lemma 3.3.}
For $h\in\Sigma_{1}$, $J_{0,0}^{'}(h)\neq0$, $J_{0,1}^{'}(h)\neq0$.
Let
$$
\omega(h)=\frac{J_{1,0}^{'}(h)}{J_{0,0}^{'}(h)},~~~\nu(h)=\frac{Z^{'}(h)}{J_{0,1}^{'}(h)},
$$
then $\omega(h)$ and $\nu(h)$ satisfy the following Riccati functions respectively
$$
D_{1}(h)\omega^{'}(h)=-g_{12}\omega^{2}(h)+(g_{22}-g_{11})\omega(h)+g_{21},
\eqno(3.21)$$
$$
D_{2}(h)\nu^{'}(h)=-k_{12}\nu^{2}(h)+(k_{32}-k_{11})\nu(h)+k_{31}.
\eqno(3.22)
$$

\noindent{\bf Proof.}
By direct computation, the trajectory passing through two saddle points corresponding to $h=0$. Suppose $x_{1}$ and $x_{2}$  are solutions of $H(x,0)=0~(x_{1}>x_{2})$. Then for $h\in\Sigma_{1}$, $x_{A}\in(0,1)$, $x_{B}\in(1,x_{1})$, and we have
$$
J_{0,0}(h)=\int_{L_{h}^{+}}dx=\int_{x_{A}(h)}^{x_{B}(h)}dx=x_{B}(h)-x_{A}(h),
$$
$$
J_{0,0}^{'}(h)=\frac{\partial{x_{B}(h)}}{\partial{h}}-\frac{\partial{x_{A}(h)}}{\partial{h}}=\frac{(x_{A}(h)-1)(x_{A}(h)-\frac{\lambda-2}{\lambda})-(x_{B}(h)-1)(x_{B}(h)-\frac{\lambda-2}{\lambda})}{(x_{B}(h)-1)(x_{B}(h)-\frac{\lambda-2}{\lambda})(x_{A}(h)-1)(x_{A}(h)-\frac{\lambda-2}{\lambda})},
$$
because of $(x_{A}(h)-1)(x_{A}(h)-(\lambda-2)/\lambda)<0$, $(x_{B}(h)-1)(x_{B}(h)-(\lambda-2)/\lambda)>0$, we can obtain $J_{0,0}^{'}(h)\neq0$.
By direct computation, we have
$$
J_{0,1}^{'}(h)=\int_{L_{h}^{+}}\frac{1}{2xy}dx=\int_{L_{h}^{+}}dt\neq0
$$
in $h\in\Sigma_{1}$.
Because of
$$
\omega^{'}(h)=\frac{1}{(J_{0,0}^{'}(h))^{2}}\left(J_{1,0}^{''}(h)J_{0,0}^{'}(h)-J_{1,0}^{'}(h)J_{0,0}^{''}(h)\right),
$$
$$
\nu^{'}(h)=\frac{1}{(J_{0,1}^{'}(h))^{2}}\left(Z^{''}(h)J_{0,1}^{'}(h)-Z^{'}(h)J_{0,1}^{''}(h)\right),
$$
by (3.16) and (3.17) we finish the proof. $\diamondsuit$

\vskip 0.3 true cm
Substituting (3.18) into (3.4), we have
$$
M(h)=\alpha(h)J_{0,0}(h)+\beta(h)J_{1,0}(h)+\gamma(h)J_{0,2}(h)+\eta(h)J_{0,1}(h)+\xi(h)J_{1,1}(h)+\zeta(h)Z(h),
$$
where the polynomial coefficients of generating elements are different from (3.4), we still use $\alpha(h)-\zeta(h)$ to express.
By (3.12) and (3.13), we have
$$
\begin{aligned}
M(h)=\alpha_{1}J_{0,0}^{'}(h)+\beta_{1}J_{1,0}^{'}(h)+\gamma_{1}J_{0,2}^{'}(h)+\eta_{1}J_{0,1}^{'}(h)+\xi_{1}J_{1,1}^{'}(h)+\zeta_{1}J_{2,1}^{'}(h),\\
M^{'}(h)=\alpha_{2}J_{0,0}^{'}(h)+\beta_{2}J_{1,0}^{'}(h)+\gamma_{2}J_{0,2}^{'}(h)+\eta_{2}J_{0,1}^{'}(h)+\xi_{2}J_{1,1}^{'}(h)+\zeta_{2}J_{2,1}^{'}(h),
\end{aligned}
\eqno(3.23)$$
where $\alpha_{s}-\zeta_{s}$ are polynomials of $h$ with ${\rm deg}\alpha_{s}(h)\leq\left[\frac{n}{3}\right]+2-s, ~{\rm deg}\beta_{s}(h)\leq\left[\frac{n-1}{3}\right]+2-s, ~{\rm deg}\gamma_{s}(h)\leq\left[\frac{n-2}{3}\right]+2-s, ~{\rm deg}\eta_{s}(h)\leq\left[\frac{n-1}{3}\right]+2-s, ~{\rm deg}\xi_{s}(h)\leq\left[\frac{n-2}{3}\right]+2-s, ~{\rm deg}\zeta_{s}(h)\leq\left[\frac{n-3}{3}\right]+2-s~(s=1,2).
$ Removing $J_{0,2}^{'}(h)$ from (3.23) gives
$$
\gamma_{1}(h)M(h)=\gamma_{2}(h)M^{'}(h)+F_{1}(h),
$$
where
$$
F_{1}(h)=\alpha_{3}(h)J_{0,0}^{'}(h)+\beta_{3}(h)J_{1,0}^{'}(h)+\eta_{3}(h)J_{0,1}^{'}(h)+\xi_{3}(h)J_{1,1}^{'}(h)+\zeta_{3}(h)J_{2,1}^{'}(h),
$$
with ${\rm deg}\alpha_{3}(h)\leq\left[\frac{n}{3}\right]+\left[\frac{n-2}{3}\right]+1,~{\rm deg}\beta_{3}(h)\leq\left[\frac{n-1}{3}\right]+\left[\frac{n-2}{3}\right]+1,~{\rm deg}\eta_{3}(h)\leq\left[\frac{n-1}{3}\right]+\left[\frac{n-2}{3}\right]+1,~{\rm deg}\xi_{3}(h)\leq2\left[\frac{n-2}{3}\right],~{\rm deg}\zeta_{3}(h)\leq\left[\frac{n-3}{3}\right]+\left[\frac{n-2}{3}\right]+1$.
By lemma $5.1$ of [15], we have
$$
\begin{aligned}
\#\{M(h)=0,h\in\Sigma_{1}\}&\leq \#\{\gamma_{1}(h)=0,h\in\Sigma_{1}\}\\
&+\#\{F_{1}(h)=0,h\in\Sigma_{1}\}+1.
\end{aligned}
\eqno(3.24)$$

\vskip 0.3 true cm
\noindent{\bf Lemma 3.4.}
Suppose that $D_{1}(h)\neq0$ for $h\in\Sigma_{1}$. Let
$$
\Phi_{1}(h)=\alpha_{3}(h)J_{0,0}^{'}(h)+\beta_{3}(h)J_{1,0}^{'}(h),
$$
$$
\Phi_{2}(h)=\eta_{3}(h)J_{0,1}^{'}(h)+\xi_{3}(h)J_{1,1}^{'}(h)+\zeta_{3}(h)Z^{'}(h).
$$
For $h\in\Sigma_{1}$, there exists polynomials $P_{2}(h)$, $P_{1}(h)$ and $P_{0}(h)$ of $h$ satisfy
${\rm deg}P_{2}(h)\leq m_{2}$, ${\rm deg}P_{1}(h)\leq m_{1}$,
${\rm deg}P_{0}(h)\leq m_{0}$
such that $$L(h)\Phi_{1}(h)=0,$$
where
$m_{2}=\left[\frac{n}{3}\right]+\left[\frac{n-1}{3}\right]+2\left[\frac{n-2}{3}\right]+13$, $m_{1}=m_{2}-1$, $m_{0}=m_{2}-2$, and
$$
L(h)=P_{2}(h)\frac{d^{2}}{dh^{2}}+P_{1}(h)\frac{d}{dh}+P_{0}(h).
\eqno(3.25)$$

\noindent{\bf Proof.} We first assert that
$$\left\{\begin{array}{c}
\Phi_{1}^{''}(h)=\frac{1}{D_{1}^{2}(h)}\left(Q_{\left[\frac{n}{3}\right]+\left[\frac{n-2}{3}\right]+5}(h)J_{0,0}^{'}(h)+Q_{\left[\frac{n-1}{3}\right]+\left[\frac{n-2}{3}\right]+5}J_{1,0}^{'}(h)\right),\\
\Phi_{1}^{'}(h)=\frac{1}{D_{1}(h)}\left(Q_{\left[\frac{n}{3}\right]+\left[\frac{n-2}{3}\right]+3}J_{0,0}^{'}(h)+Q_{\left[\frac{n-1}{3}\right]+\left[\frac{n-2}{3}\right]+3}(h)J_{1,0}^{'}(h)\right),
\end{array}\right.
$$
where $Q_{s}(h)$ express ${\rm deg}Q_{s}(h)\leq s$. In fact, by (3.16), we have
$$
\begin{aligned}
\Phi_{1}^{'}(h)&=\alpha_{3}^{'}(h)J_{0,0}^{'}(h)+\alpha_{3}(h)J_{0,0}^{''}(h)+\beta_{3}^{'}(h)J_{1,0}^{'}(h)+\beta_{3}(h)J_{1,0}^{''}(h)\\
&=\left(\alpha_{3}^{'}(h)+\frac{1}{D_{1}(h)}\alpha_{3}(h)g_{11}(h)+\frac{1}{D_{1}(h)}\beta_{3}(h)g_{21}(h)\right)J_{0,0}^{'}(h)\\
&+\left(\beta_{3}^{'}(h)+\frac{1}{D_{1}(h)}\alpha_{3}(h)g_{12}(h)+\frac{1}{D_{1}(h)}\beta_{3}(h)g_{22}(h)\right)J_{1,0}^{'}(h)\\
&:=\frac{1}{D_{1}(h)}\left(Q_{\left[\frac{n}{3}\right]+\left[\frac{n-2}{3}\right]+3}J_{0,0}^{'}(h)+Q_{\left[\frac{n-1}{3}\right]+\left[\frac{n-2}{3}\right]+3}(h)J_{1,0}^{'}(h)\right).
\end{aligned}
$$
The result for $\Phi_{1}^{''}(h)$ can be proved similarly.
Suppose that
$$
P_{2}(h)=\sum_{k=0}^{m_{2}}p_{2,k}h^{k},~~P_{1}(h)=\sum_{k=0}^{m_{1}}p_{1,k}h^{k},~~P_{0}(h)=\sum_{k=0}^{m_{0}}p_{0,k}h^{k},
$$
are polynomials of $h$ with coefficients $p_{2,k}$, $p_{1,k}$ and $p_{0,k}$.
Then by the process of simplification, we have
$$
L(h)\Phi_{1}(h)=\frac{1}{D_{1}^{2}(h)}\left(X(h)J_{0,0}^{'}(h)+Y(h)J_{1,0}^{'}(h)\right),
$$
where $X(h)$ and $Y(h)$ are polynomials of $h$ with degree no more then $2\left[\frac{n}{3}\right]+\left[\frac{n-1}{3}\right]+3\left[\frac{n-2}{3}\right]+18$ and $\left[\frac{n}{3}\right]+2\left[\frac{n-1}{3}\right]+3\left[\frac{n-2}{3}\right]+18$ respectively. Let
$$
X(h)=\sum_{i=0}^{m_{3}}x_{i}h^{i},~~~Y(h)=\sum_{j=0}^{m_{4}}y_{j}h^{j},
$$
$x_{i}$ and $y_{j}$ are expressed by $p_{2,k}$, $p_{1,k}$ and $p_{0,k}$ linearly. So $L(h)F_{1}(h)=0$ is satisfied if we let
$$
x_{i}=0,~~~y_{j}=0,~~~(0\leq i\leq m_{3},~0\leq j\leq m_{4}).
\eqno(3.26)$$
System (3.26) is a homogeneous linear equation with $3\left[\frac{n}{3}\right]+3\left[\frac{n-1}{3}\right]+6\left[\frac{n-2}{3}\right]+38$ equations about $3\left[\frac{n}{3}\right]+3\left[\frac{n-1}{3}\right]+6\left[\frac{n-2}{3}\right]+39$ variables of $p_{2,k}$, $p_{1,k}$ and $p_{0,k}$. Since it follows that from the theory of linear algebra that there exist $p_{2,k}$, $p_{1,k}$ and $p_{0,k}$ such that (3.26) holds, which yields the desired result. $\diamondsuit$

\vskip 0.1 true cm
By (3.17), we have
$$\left\{\begin{array}{c}
\Phi_{2}^{''}(h)=\frac{1}{D_{2}^{2}(h)}\left(Q_{\left[\frac{n-1}{3}\right]+\left[\frac{n-2}{3}\right]+5}(h)J_{0,1}^{'}(h)+Q_{\left[\frac{n-3}{3}\right]
+\left[\frac{n-2}{3}\right]+5}(h)Z^{'}(h)+Q_{2\left[\frac{n-2}{3}\right]+5}(h)J_{1,1}^{'}(h)\right),\\
\Phi_{2}^{'}(h)=\frac{1}{D_{2}(h)}\left(Q_{\left[\frac{n-1}{3}\right]+\left[\frac{n-2}{3}\right]+3}(h)J_{0,1}^{'}(h)+Q_{\left[\frac{n-3}{3}\right]
+\left[\frac{n-2}{3}\right]+3}(h)Z^{'}(h)+Q_{2\left[\frac{n-2}{3}\right]+3}(h)J_{1,1}^{'}(h)\right),
\end{array}\right.
$$
where $Q_{s}(h)$ express ${\rm deg}Q_{s}(h)\leq s$. Then it is similar to the proof of Lemma 3.4, we have the following result.

\vskip 0.2 true cm
\noindent{\bf Lemma 3.5.}
Suppose that $D_{2}(h)\neq0$ for $h\in\Sigma_{1}$, then we have $L(h)F_{1}(h)=R(h)$, where $L(h)$ is defined as (3.25), and
$$
R(h)=\frac{1}{D_{2}^{2}(h)}\left(Q_{1}(h)J_{0,1}^{'}(h)+Q_{2}(h)Z^{'}(h)+Q_{3}(h)J_{1,1}^{'}(h)\right),
\eqno(3.27)$$
$Q_{s}(h)$ $(s=1,2,3)$ are polynomials of $h$ with degree no more than
$\left[\frac{n}{3}\right]+2\left[\frac{n-1}{3}\right]+3\left[\frac{n-2}{3}\right]+18$,
$\left[\frac{n}{3}\right]+3\left[\frac{n-2}{3}\right]+\left[\frac{n-1}{3}\right]+\left[\frac{n-3}{3}\right]+18$,
$\left[\frac{n}{3}\right]+\left[\frac{n-1}{3}\right]+4\left[\frac{n-2}{3}\right]+18$ respectively.

\vskip 0.3 true cm
\noindent{\bf Lemma 3.6.}
For $h\in\Sigma_{1}$, the number of zeros of $\Phi_{1}(h)$ does not exceed $\left[\frac{n}{3}\right]+2\left[\frac{n-1}{3}\right]+3\left[\frac{n-2}{3}\right]+9$; the number of zeros of $R(h)$ does not exceed $7\left[\frac{n}{3}\right]+11\left[\frac{n-1}{3}\right]+25\left[\frac{n-2}{3}\right]+2\left[\frac{n-3}{3}\right]+141$ (counting the multiplicity).

\vskip 0.1 true cm
\noindent{\bf Proof.} Let
$$
S_{1}(h)=\frac{\Phi_{1}(h)}{J_{0,0}^{'}(h)}=\alpha_{3}(h)+\beta_{3}(h)\omega(h).
$$
By Lemma $3.3$, we have
$$
\beta_{3}(h)D_{1}(h)S_{1}^{'}(h)=-g_{12}(h)S_{1}^{2}(h)+N_{1}(h)S_{1}(h)+N_{2}(h),
$$
where
$$
N_{1}(h)=D_{1}(h)\beta_{3}^{'}(h)+2\alpha_{3}(h)g_{12}(h)+\beta_{3}(h)(g_{22}(h)-g_{11}(h)),
$$
$$
\begin{aligned}
N_{2}(h)&=\alpha_{3}(h)\left(D_{1}(h)\beta_{3}^{'}(h)+2\alpha_{3}(h)g_{12}(h)+\beta_{3}(h)\left(g_{22}(h)-g_{11}(h)\right)\right)+g_{12}(h)\alpha_{3}^{2}(h)\\
&+\beta_{3}^{2}(h)g_{21}(h)+\beta_{3}(h)D_{1}(h)\beta_{3}^{'}(h).
\end{aligned}
$$
Then by Lemma 5.1 of [20], we can obtain
$$
\begin{aligned}
\#\{\Phi_{1}(h)=0,h\in\Sigma_{1}\}&=\#\{S_{1}(h)=0,h\in\Sigma_{1}\}\leq\#\{\beta_{3}(h)D_{1}(h)=0,h\in\Sigma_{1}\}\\
&+\#\{N_{2}(h)=0,h\in\Sigma_{1}\}+1\leq \left[\frac{n}{3}\right]+2\left[\frac{n-1}{3}\right]+3\left[\frac{n-2}{3}\right]+9.
\end{aligned}
\eqno(3.28)$$
Let
$$
R_{1}(h)=Q_{1}(h)J_{0,1}^{'}(h)+Q_{2}(h)Z^{'}(h)+Q_{3}(h)J_{1,1}^{'}(h),
$$
suppose $\Delta$ is the set of zeros of $Q_{3}(h)$ in $(-2,0)$, then in $(-2,0)\backslash\ \Delta$, we can get
$$
\frac{R_{1}(h)}{Q_{3}(h)}=\frac{1}{Q_{3}(h)}\left(Q_{1}(h)J_{0,1}^{'}(h)+Q_{2}(h)Z^{'}(h)\right)+J_{1,1}^{'}(h).
$$
Using (3.17), we get that
$$
\left(\frac{R_{1}(h)}{Q_{3}(h)}\right)^{'}=\frac{F_{2}(h)}{D_{2}(h)Q_{3}^{2}(h)},
$$
where
$$
F_{2}(h)=\eta_{4}(h)J_{0,1}^{'}(h)+\zeta_{4}(h)Z^{'}(h),
$$
with ${\rm deg}\eta_{4}(h)\leq2\left[\frac{n}{3}\right]+3\left[\frac{n-1}{3}\right]+7\left[\frac{n-2}{3}\right]+38$, ${\rm deg}\zeta_{4\frac{}{}}(h)\leq2\left[\frac{n}{3}\right]+2\left[\frac{n-1}{3}\right]+7\left[\frac{n-2}{3}\right]+\left[\frac{n-3}{3}\right]+38$.
By Lemma 4.1 and Lemma 4.2 of [27], we can obtain
$$
\#\{R_{1}(h)=0,h\in\Sigma_{1}\}\leq\#\{F_{2}(h)=0,h\in\Sigma_{1}\}+\#\{Q_{3}(h)=0,h\in\Sigma_{1}\}+3.
$$
Let
$$
S_{2}(h)=\frac{F_{2}(h)}{J_{0,1}^{'}(h)}=\eta_{4}(h)+\zeta_{4}(h)\nu(h),
$$
similar to the upper bound of zeros number of the $\Phi_{1}(h)$, using Lemma 3.3, we get that
$$
\begin{aligned}
\#\{F_{2}(h)=0,h\in\Sigma_{1}\}&=\#\{S_{2}(h)=0,h\in\Sigma_{1}\}\\
&\leq6\left[\frac{n}{3}\right]+7\left[\frac{n-1}{3}\right]+21\left[\frac{n-2}{3}\right]+2\left[\frac{n-3}{3}\right]+120.
\end{aligned}
$$
Therefore, we have
$$
\begin{aligned}
\#\{R(h)=0,h\in\Sigma_{1}\}&=\#\{R_{1}(h)=0,h\in\Sigma_{1}\}\\
&\leq7\left[\frac{n}{3}\right]+8\left[\frac{n-1}{3}\right]+25\left[\frac{n-2}{3}\right]+2\left[\frac{n-3}{3}\right]+141.
\end{aligned}
\eqno(3.29)$$
This ends the proof. $\diamondsuit$

\vskip 0.5 true cm
\noindent{\bf Proof for the elliptic segment.}
\vskip 0.3 true cm
We denote $m_{5}=\left[\frac{n}{3}\right]+2\left[\frac{n-1}{3}\right]+3\left[\frac{n-2}{3}\right]+9$, $m_{6}=7\left[\frac{n}{3}\right]+11\left[\frac{n-1}{3}\right]+25\left[\frac{n-2}{3}\right]+2\left[\frac{n-3}{3}\right]+141$. For $h\in\Sigma_{1}$, by Lemma 3.4 and Lemma 3.6 we can assume that
$$
P_{2}(\widetilde{h}_{i})=0,~~\Phi_{1}(\overline{h}_{j})=0,~~\widetilde{h}_{i},\overline{h}_{j}\in\Sigma_{1},~~1\leq i\leq m_{2} ,~~1\leq j\leq m_{5}.
$$
Denote $\widetilde{h}_{i}$ and $\overline{h}_{j}$ as $h_{k}^{*}<h_{k+1}^{*}$ for $k=1,2,...,m_{2}+m_{5}$. Let
$$
\Delta_{s}=(h_{s}^{*},h_{s+1}^{*}),~~s=0,1,...,m_{2}+m_{5},
$$
where $h_{0}^{*}=\lambda-3$ and $h_{m_{2}+m_{5}+1}^{*}=0$. Then $P_{2}(h)\neq0$ and $\Phi_{1}(h)\neq0$ for $h\in\Delta_{s}$ and $L(h)\Phi_{1}(h)=0$. By Lemma 2.3, the solution space of
$$
L(h)=P_{2}(h)\left(\frac{d^{2}}{dh^{2}}+\frac{P_{1}(h)}{P_{2}(h)}\frac{d}{dh}+\frac{P_{0}(h)}{P_{2}(h)}\right)
$$
is a Chebyshev space on $\Delta_{s}$. By Lemma 2.4, $M(h)$ has at most $2+l_{s}$ zeros for $h\in\Delta_{s}$, where $l_{s}$ is the number of zeros of $R(h)$ on $\Delta_{s}$. Then we obtain
$$
\begin{aligned}
\#\{M(h)=0,h\in\Sigma_{1}\}&\leq \#\{\gamma_{1}(h)=0,h\in\Sigma_{1}\}+\#\{F_{1}(h)=0,h\in\Sigma_{1}\}+1\\
&\leq \#\{R(h)=0,h\in\Sigma_{1}\}+2\cdot the ~number ~of ~the ~intervals ~of ~\Delta_{s}\\
&+the ~number ~of ~the ~end ~points ~of ~\Delta_{s}+\#\{\gamma_{1}(h)=0,h\in\Sigma_{1}\}+1\\
&\leq m_{6}+3(m_{2}+m_{5})+\left[\frac{n-2}{3}\right]+4\\
&\leq 25n+161.
\end{aligned}
$$
By the same arguments as the proof as above, we get
$$
\#\{M(h)=0,h\in\Sigma_{2}\}\leq 25n+161.
$$
This ends the proof.

\vskip 0.3 true cm
\noindent\quad{{\bf $\S$3.2 The case of hyperbolic segment.}}
\vskip 0.3 true cm
In is case, we can get (3.1) and (3.2) by the same transformation of coordinates of elliptic segment.

For $\lambda\in(-1,0)$, system (3.2) has an elementary center (1,0) corresponding to $h=\lambda-3$, three saddles $(0,\pm\sqrt{3(\lambda-2)})$ and $((\lambda-2)/\lambda,0)$ corresponding to $h=0$ and $h=(\lambda-2)^{2}(\lambda+1)/\lambda^{2}$ respectively, then there exist period annulus if $h\in(\lambda-3,0)$ (see Fig.\,2(b)).

By the same arguments as the proof as the case of elliptic segment, we can obtain the first order Melnikov function $M(h)$ of this case satisfy
$$
\#\{M(h)=0,h\in(\lambda-3,0)\}\leq 25n+161.$$

\vskip 0.3 true cm
\noindent\quad{{\bf $\S$3.3 The case of parabolic segment.}}
\vskip 0.3 true cm
In this case, let
$$x_{1}=2\sqrt{2}y-2x+2,~~y_{1}=\sqrt{2}x+y,~~dt_{1}=\frac{1}{2}dt.$$
Then (1.1) and (1.2) transform into
$$H(x,y)=\frac{1}{2}x(2y^2+x-4)=h, ~~~h\in(-2,0),\eqno(3.30)$$
$$
\left(
  \begin{array}{c}
          \dot{x}\\
          \dot{y}
   \end{array}
   \right)=\begin{cases}
   \left(
    \begin{array}{c}
        -2xy+\varepsilon p^+(x,y)\\
        -2+x+y^2+\varepsilon q^+(x,y)
   \end{array}
   \right), \quad y>0,\\
   \left(
    \begin{array}{c}
       -2xy+\varepsilon p^-(x,y)\\
      -2+x+y^2+\varepsilon q^-(x,y)
   \end{array}
   \right), \quad y<0.\\
   \end{cases}
   \eqno(3.31)$$
here we shall omit the subscript 1. The point $(2,0)$ is an elementary center corresponding to $h=-2$ and the point $(0,-\sqrt{2})$, $(0,\sqrt{2})$ are saddles corresponding to $h=0$, then there exist period annulus if $h\in(-2,0)$ (see Fig.\,2(c)).

\vskip 0.3 true cm
\noindent{\bf Lemma 3.7.}
If $n\geq2$, then for $h\in(-2,0)$, $M(h)$ can be expressed as
$$
M(h)=\alpha(h)J_{0,1}(h)+\beta(h)J_{1,1}(h)+\gamma(h)J_{1,0}(h)+\eta(h)J_{0,2}(h),
\eqno(3.32)
$$
where
$$
{\rm deg}\alpha(h)\leq\left[\frac{n}{2}\right],~~{\rm deg}\beta(h)\leq\left[\frac{n}{2}\right]-1,
~~{\rm deg}\gamma(h)\leq\left[\frac{n}{2}\right],~~{\rm deg}\eta(h)\leq\left[\frac{n-2}{3}\right].
$$

\vskip 0.1 true cm
\noindent{\bf Proof.}
Similar to the proof of Lemma 3.1, we can get
$$
M(h)=\sum\limits_{\begin{subarray}{c}i+j=0\\i\geq0,j\geq0\end{subarray}}^{n}\rho_{i,j}J_{i,j}(h),
$$
$$
L_{h}^{+}(L_{h}^{-})=\{(x,y)|H(x,y)=h,y>0(y<0),h\in(-2,0)\}.
$$
Differentiating (3.30) with respect to $x$, we obtain
$$
y^{2}+2xy\frac{\partial{y}}{\partial{x}}+x-2=0.\eqno(3.33)
$$
Multiplying $(3.33)$ by $x^{i-1}y^{j}dx$, integrating over $L_{h}^{+}$, we have
$$
\frac{j-2i-2}{j}J_{i,j}(h)=2J_{i,j-2}(h)-J_{i+1,j-2}(h).
\eqno(3.34)$$
Multiplying $(3.30)$ by $x^{i}y^{j}dx$, integrating over $L_{h}^{+}$, we have
$$
J_{i,j}(h)=hJ_{i-1,j-2}(h)+2J_{i,j-2}(h)-\frac{1}{2}J_{i+1,j-2}(h).
\eqno(3.35)$$
Noting that we have
$$
2J_{i,2i}(h)=J_{i+1,2i}(h),
$$
when $j-2i-2=0$ in $(3.34)$.
By $(3.34)$, we have
$$
M(h)=\sum\limits_{i=0}^{n-1}c_{i,1}J_{i,1}(h)+\sum\limits_{i=0}^{n}c_{i,0}J_{i,0}(h)+\sum\limits_{i=0}^{\left[\frac{n-2}{3}\right]}c_{i,2i+2}J_{i,2i+2}(h),
\eqno(3.36)$$
where $c_{i,1}$, $c_{i,0}$, $c_{i,2i+2}$ are some real constants.
By $(3.34)$ and $(3.35)$, we have
$$
h(j-2i)J_{i-2,j-2}+2(j-2i)J_{i-1,j-2}-\frac{1}{2}(j-2i)J_{i,j-2}=2jJ_{i-1,j-2}-jJ_{i,j-2}.
\eqno(3.37)$$
Let $j=2$, $j=3$ in $(3.37)$, we have respectively
$$
\begin{aligned}
(i+1)J_{i,0}=4iJ_{i-1,0}+2(i-1)hJ_{i-2,0},\\
(\frac{3}{2}+i)J_{i,0}=4iJ_{i-1,1}-h(3-2i)J_{i-2,1}.
\end{aligned}
\eqno(3.38)$$
Let $j=2i+2$ in $(3.35)$, we have
$$
J_{i,2i+2}=hJ_{i-1,2i}+2J_{i,2i}-\frac{1}{2}J_{i+1,2i}.
\eqno(3.39)$$
By $(3.34)-(3.39)$, we have
$$
M(h)=\alpha(h)J_{0,1}(h)+\beta(h)J_{1,1}(h)+\gamma(h)J_{1,0}(h)+\eta(h)J_{0,2}(h),
$$
where
$$
{\rm deg}\alpha(h)\leq\left[\frac{n}{2}\right],~~{\rm deg}\beta(h)\leq\left[\frac{n}{2}\right]-1,
{\rm deg}\alpha(h)\leq\left[\frac{n}{2}\right],~~{\rm deg}\eta(h)\leq\left[\frac{n-2}{3}\right].
$$
This ends the proof. $\diamondsuit$

Similar to the proof of Lemma 3.2, Lemma 3.4 and Lemma 3.5, we can get the following Lemmas.

\vskip 0.3 true cm
\noindent{\bf Lemma 3.8.}
Let ${\bf U_{1}(h)}=(J_{0,1}(h),J_{1,1}(h))^{T}$, ${\bf U_{2}(h)}=(J_{1,0}(h),J_{0,2}(h))^{T}$,
than the vector functions ${\bf U_{1}(h)}$ and ${\bf U_{2}(h)}$ satisfy respectively the following Picard-Fuchs equations:
$$
{\bf U_{1}(h)}=(B_{1}h+C_{1}){\bf U_{1}'(h)},\eqno(3.40)
$$
$$
{\bf U_{2}(h)}=(B_{2}h+C_{2}){\bf U_{2}'(h)},\eqno(3.41)
$$
where
\begin{eqnarray*}
B_{1}h+C_{1}=\left(\begin{matrix}
                         \frac{4}{3}h&\frac{4}{3}\\
                         \frac{8}{15}h&\frac{4}{5}h+\frac{32}{15}
                         \end{matrix}\right),
~~~
B_{2}h+C_{2}=\left(\begin{matrix}
                         2h+4&0\\
                         h+2&h
                         \end{matrix}\right).
\end{eqnarray*}
And we have
$$
D_{1}(h):={\rm det}(B_{1}h+C_{1})=\frac{16}{15}h(h+2),
$$
$$
D_{2}(h):={\rm det}(B_{2}h+C_{2})=2h(h+2).
$$

\vskip 0.3 true cm
\noindent{\bf Lemma 3.9.}
Let
$$
\Phi_{1}(h)=\alpha(h)J_{0,1}(h)+\beta(h)J_{1,1}(h),~~~\Phi_{2}(h)=\gamma(h)J_{1,0}(h)+\eta(h)J_{0,2}(h).
$$
For $h\in(-2,0)$, there exist polynomials $P_{2}(h)$, $P_{1}(h)$ and $P_{0}(h)$ of $h$ with degree respectively $2\left[\frac{n}{2}\right]+2$, $2\left[\frac{n}{2}\right]+1$ and $2\left[\frac{n}{2}\right]$ such that $L(h)\Phi_{1}(h)=0$, where
$$
L(h)=P_{2}(h)\frac{d^{2}}{dh^{2}}+P_{1}(h)\frac{d}{dh}+P_{0}(h).
\eqno(3.42)$$

\vskip 0.1 true cm
By (3.41) we have
$$
-J_{1,0}^{'}(h)=(4+2h)J_{1,0}^{''}(h),~~~-J_{1,0}^{'}(h)=(h+1)J_{1,0}^{''}(h)+hJ_{0,2}^{''}(h).
\eqno(3.43)$$
Hance we have
$$hJ_{0,2}^{''}(h)=(3+h)J_{1,0}^{''}(h).
\eqno(3.44)$$
Noticing (3.43) and (3.44), we can get the following result as similar with Lemma 3.5.
\vskip 0.3 true cm
\noindent{\bf Lemma 3.10.} Suppose that $D_{2}(h)\neq0$ for $h\in(-2,0)$, then we have $L(h)M(h)=R(h)$, where $L(h)$ is defined as (3.42), and
$$
R(h)=\frac{1}{h}\left(M_{1}(h)J_{1,0}^{''}(h)+M_{2}(h)J_{0,2}^{'}(h)\right),
\eqno(3.45)$$
$M_{1}(h)$, $M_{2}(h)$ are polynomials of $h$ with degree no more than $3\left[\frac{n}{2}\right]+3$, $2\left[\frac{n}{2}\right]+\left[\frac{n-2}{3}\right]+2$ respectively.

\vskip 0.3 true cm
\noindent{\bf Lemma 3.11.}
When $h\in(-2,0)$, $\Phi_{1}(h)$, $R(h)$ has at most $3\left[\frac{n}{2}\right]+2$, $6\left[\frac{n}{2}\right]+3\left[\frac{n-2}{3}\right]+12$ zeros respectively (taking into account the multiplicity).

\vskip 0.1 true cm
\noindent{\bf Proof.}
Let
$$
\sigma(h)=(\alpha(h),\beta(h)),~~~\tau(h)=(\gamma(h),\eta(h)).
$$
By $(3.40)$, we have
$$
{\rm det}(B_{1}h+C_{1})U_{1}^{'}(h)=(B_{1}h+C_{1})^{*}U_{1}(h),
\eqno(3.46)$$
where $(B_{1}h+C_{1})^{*}:=(g_{ij})_{2\times2}$.
Notice that
$$
J_{0,1}(h)=\int_{L_{h}^{+}}ydx=\int_{L_{h}^{+}}ydx+\int_{\overrightarrow{B_{h}A_{h}}}ydx=\iint_{D}dxdy=S_{D}\neq0,
$$
where $S_{D}$ express the area of $D$. Let $\omega(h)=\frac{J_{1,1}(h)}{J_{0,1}(h)}$, by $(3.46)$ we can obtain $\omega(h)$ satisfy the following Riccati equation
$$
D_{1}(h)\omega^{'}(h)=-g_{12}\omega^{2}(h)+(g_{22}-g_{11})\omega(h)+g_{21}.
$$
Let $S_{1}(h)=\frac{\Phi_{1}(h)}{J_{0,1}(h)}=\alpha(h)+\beta(h)\omega(h)$, then $S_{1}(h)$ satisfy the following Riccati equation
$$
\beta(h)D_{1}(h)S_{1}^{'}(h)=-g_{12}S_{1}^{2}(h)+N_{1}(h)S_{1}(h)+N_{2}(h),
$$
where $N_{1}(h)$ and $N_{2}(h)$ are polynomials of $h$ with degree no more than $\left[\frac{n}{2}\right]$, $2\left[\frac{n}{2}\right]$ respectively.
Then using Lemma 5.1 of [20], we have
$$
\begin{aligned}
\#\{\Phi_{1}(h)=0,h\in(-2,0)\}&=\#\{S_{1}(h)=0,h\in(-2,0)\}\\
&\leq\#\{\beta(h)D_{1}(h)=0,h\in(-2,0)\}+\#\{N_{2}(h)=0,h\in(-2,0)\}+1\\
&\leq3\left[\frac{n}{2}\right]+2.
\end{aligned}
$$
Let
$$
R_{1}(h)=hR(h)=M_{1}(h)J_{1,0}^{''}(h)+M_{2}(h)J_{0,2}^{'}(h),
$$
suppose $\Delta$ is the set of zeros of $M_{2}(h)$ in $(-2,0)$, then in $(-2,0)\backslash\Delta$, we can get
$$
\left(\frac{R_{1}(h)}{M_{2}(h)}\right)^{'}=\frac{F_{3}(h)}{h(4+2h)M_{2}^{2}(h)},~~F_{3}(h)=M_{3}(h)J_{1,0}^{''}(h),
$$
where $M_{3}(h)$ is the polynomials of $h$ with degree no more than $4\left[\frac{n}{2}\right]+2\left[\frac{n-2}{3}\right]+6$ in $h\in(-2,0)$. By Lemma 4.1 and Lemma 4.2 of [27], we have
$$
\#\{R_{1}(h)=0,h\in(-2,0)\}\leq\#\{F_{3}(h)=0,h\in(-2.0)\}+\#\{M_{2}(h)=0,h\in(-2,0)\}+3.
$$
By direct computation, we can obtain $J_{1,0}^{''}(h)=-\frac{C}{4(h+2)^{\frac{3}{2}}}$ ($C$ is a constant), then
$$
F_{3}(h)=M_{3}(h)\frac{C}{4(h+2)}\frac{1}{\sqrt{h+2}}.
$$
Using Lemma 2.5 yields
$$
\#\{F_{3}(h)=0,h\in(-2,0)\}\leq\#\{M_{3}(h)=0,h\in(-2,0)\}+1\leq4\left[\frac{n}{2}\right]+2\left[\frac{n-2}{3}\right]+7.
$$
Hance,
$$
\#\{R(h)=0,h\in(-2,0)\}=\#\{R_{1}(h)=0,h\in(-2,0)\}\leq6\left[\frac{n}{2}\right]+3\left[\frac{n-2}{3}\right]+12.
$$
This ends the proof. $\diamondsuit$

\vskip 0.3 true cm
\noindent{\bf Proof for the case of parabolic segment.}

For $h\in(-2,0)$, by Lemma 3.9 and Lemma 3.11 we can assume that
$$
P_{2}(\widetilde{h}_{i})=0,~~\Phi_{1}(\overline{h}_{j})=0,~~\widetilde{h}_{i},\overline{h}_{j}\in(-2,0),~~1\leq i\leq 2\left[\frac{n}{2}\right]+2 ,~~1\leq j\leq 3\left[\frac{n}{2}\right]+2.
$$
Denote $\widetilde{h}_{i}$ and $\overline{h}_{j}$ as $h_{k}^{*}<h_{k+1}^{*}$ for $k=1,2,...,5\left[\frac{n}{2}\right]+4$. Let
$$
\Delta_{s}=(h_{s}^{*},h_{s+1}^{*}),~~s=0,1,...,5\left[\frac{n}{2}\right]+4,
$$
where $h_{0}^{*}=-2$ and $h_{5\left[\frac{n}{2}\right]+5}^{*}=0$. Then $P_{2}(h)\neq0$ and $\Phi_{1}(h)\neq0$ for $h\in\Delta_{s}$ and $L(h)\Phi_{1}(h)=0$. By Lemma 2.3, the solution space of
$$
L(h)=P_{2}(h)\left(\frac{d^{2}}{dh^{2}}+\frac{P_{1}(h)}{P_{2}(h)}\frac{d}{dh}+\frac{P_{0}(h)}{P_{2}(h)}\right)
$$
is a Chebyshev space on $\Delta_{s}$. By Lemma 2.4, $M(h)$ has at most $2+l_{s}$ zeros for $h\in\Delta_{s}$, where $l_{s}$ is the number of zeros of $R(h)$ on $\Delta_{s}$. We obtain
$$
\begin{aligned}
\#\{M(h)=0,h\in(-2,0)\}&\leq \#\{R(h)=0,h\in(-2,0)\}+2\cdot the ~number ~of ~the ~intervals ~of ~\Delta_{s}\\
&\leq the ~number ~of ~the ~end ~points ~of ~\Delta_{s}\\
&\leq 6\left[\frac{n}{2}\right]+3\left[\frac{n-2}{3}\right]+12+2(5\left[\frac{n}{2}\right]+5)+5\left[\frac{n}{2}\right]+4\\
&\leq 12n+24.
\end{aligned}$$
This ends the proof.

\vskip 0.3 true cm
\centerline{\bf { $\S$4.} Proof of the results on  $S^{(3)}$}
\vskip 0.2 true cm

In this case, $a=1$, $b=0$. (1.1) and (1.2) turns to
$$
H(x,y)=\frac{1}{2}(x^{2}+y^{2})-\frac{1}{3}x^{3}+xy^{2}=h,~~~h\in(0,\frac{1}{6}).
\eqno(4.1)$$
$$
\left(
  \begin{array}{c}
          \dot{x}\\
          \dot{y}
   \end{array}
   \right)=\begin{cases}
   \left(
    \begin{array}{c}
        y+2xy+\varepsilon p^+(x,y)\\
        -x+x^{2}-y^{2}+\varepsilon q^+(x,y)
   \end{array}
   \right), \quad y>0,\\
   \left(
    \begin{array}{c}
       y+2xy+\varepsilon p^-(x,y)\\
      -x+x^{2}-y^{2}+\varepsilon q^-(x,y)
   \end{array}
   \right), \quad y<0.\\
   \end{cases}
   \eqno(4.2)$$
System (4.2) has an elementary center point $(0,0)$ corresponding to $h=0$, three saddles points $(1,0)$, $(-\frac{1}{2},\pm\frac{\sqrt{3}}{2})$ and corresponding to $h=\frac{1}{6}$, then there exist period annulus if $h\in(0,\frac{1}{6})$ (see Fig.\,2(d)).
\vskip 0.1 true cm
Similar to the proof of Lemma 3.1, we can get the following Lemmas.
\vskip 0.3 true cm
\noindent{\bf Lemma 4.1.}
Consider system (3.2), if $n\geq3$, for $h\in(0,\frac{1}{6})$, $M(h)$ can be expressed as
$$
M(h)={\bf \sigma_{1}(h)}{\bf U_{1}(h)}+{\bf \sigma_{2}(h)}{\bf U_{2}(h)},
\eqno(4.3)$$
where
$$
{\bf \sigma_{1}(h)}=(\alpha(h),\beta(h),\gamma(h)),~~{\bf \sigma_{2}(h)}=(\eta(h),\xi(h),\zeta(h)),
$$
$$
{\bf U_{1}(h)}=(J_{0,0}(h),J_{1,0}(h),J_{0,2}(h))^{T},~~{\bf U_{2}(h)}=(J_{0,1}(h),J_{1,1}(h),J_{2,1}(h))^{T},
$$
and $\alpha(h)$, $\beta(h)$, $\gamma(h)$, $\eta(h)$, $\xi(h)$ and $\zeta(h)$ are polynomials of $h$ with
$$
{\rm deg}\alpha(h)\leq\left[\frac{n}{3}\right],~~{\rm deg}\beta(h)\leq\left[\frac{n-1}{3}\right],~~{\rm deg}\gamma(h)\leq\left[\frac{n-2}{3}\right],
$$
$$
{\rm deg}\eta(h)\leq\left[\frac{n-1}{3}\right],~~{\rm deg}\xi(h)\leq\left[\frac{n-2}{3}\right],~~{\rm deg}\zeta(h)\leq\left[\frac{n-3}{3}\right].
$$

\vskip 0.3 true cm
\noindent{\bf Lemma 4.2.}
If $n\geq3$, for $h\in(0,\frac{1}{6})$.
\vskip 0.1 true cm
\noindent(1) The vector functions ${\bf U_{1}(h)}$ and ${\bf U_{2}(h)}$ satisfy respectively the following Picard-Fuchs equation:
$$
{\bf U_{1}(h)}=(B_{1}h+C_{1}){\bf U_{1}^{'}(h)},
\eqno(4.4)$$
$$
{\bf U_{2}(h)}=(B_{2}h+C_{2}){\bf U_{2}^{'}(h)},
\eqno(4.5)$$
where
$$
B_{1}h+C_{1}=\left(\begin{matrix}
                   &3h~~&-\frac{1}{2}~~&0\\
                   &\frac{3}{4}h~~&\frac{3}{2}h-\frac{3}{8}~~&0\\
                   &\frac{3}{4}h~~&-\frac{1}{2}h-\frac{1}{24}~~&h-\frac{1}{6}
                   \end{matrix}\right),
~~~
B_{2}h+C_{2}=\left(\begin{matrix}
                   &\frac{3}{2}h~~&-\frac{1}{2}~~&\frac{1}{4}\\
                   &0~~&h-\frac{1}{6}~~&0\\
                   &\frac{3}{16}h~~&\frac{1}{8}h+\frac{1}{96}~~&\frac{3}{4}h-\frac{3}{16}
                   \end{matrix}\right).
$$
And we have
$$
D_{1}(h):={\rm det}(B_{1}h+C_{1})=\frac{1}{8}h(6h-1)^{2},
$$
$$
D_{2}(h):={\rm det}(B_{2}h+C_{2})=\frac{1}{128}h(6h-1)(24h-7).
$$

\noindent(2)
$$
{\rm det}(B_{1}h+C_{1}){\bf U_{1}^{''}(h)}=\left(\begin{matrix}
                                     &e_{11}(h)~~&e_{12}(h)\\
                                     &e_{21}(h)~~&e_{22}(h)\\
                                     &e_{31}(h)~~&e_{32}(h)
                                     \end{matrix}\right)
\left(\begin{matrix}
       J_{0,0}^{'}(h)\\
       J_{1,0}^{'}(h)
       \end{matrix}\right),
\eqno(4.6)$$
$$
\frac{3}{64}h(24h-7)\left(\begin{matrix}
                         J_{0,1}^{''}(h)\\
                         J_{2,1}^{(3)}(h)
                         \end{matrix}\right)
=\left(\begin{matrix}
        &l_{11}(h)~~&l_{12}(h)\\
        &l_{21}(h)~~&l_{22}(h)
        \end{matrix}\right)
\left(\begin{matrix}
       J_{0,1}^{'}(h)\\
       J_{2,1}^{''}(h)
       \end{matrix}\right),
\eqno(4.7)$$
where
$$
\begin{aligned}
&e_{11}(h)=-\frac{1}{8}(4h-1)(6h-1)-\frac{3}{8}h+\frac{1}{16},~~~e_{12}(h)=-\frac{1}{4}h+\frac{1}{24},\\
&e_{21}(h)=-\frac{1}{8}h(6h-1),~~~~~~~~~~~~~~~~~~~~~~~~~~~e_{22}(h)=-\frac{1}{4}h(6h-1),\\
&e_{31}(h)=-\frac{3}{2}h^{2}+\frac{1}{4}h,~~~~~~~~~~~~~~~~~~~~~~~~~~~~~e_{32}(h)=\frac{3}{2}h^{2}-\frac{1}{4}h,\\
&l_{11}(h)=-\frac{1}{64}(24h-7),~~~~~~~~~~~~~~~~~~~~~~~~~~~l_{12}(h)=-\frac{1}{128}(24h-7),\\
&l_{21}(h)=-\frac{1}{16},~~~~~~~~~~~~~~~~~~~~~~~~~~~~~~~~~~~~~~~l_{22}(h)=-\frac{3}{4}h+\frac{1}{32}.
\end{aligned}
$$

\vskip 0.1 true cm
\noindent{\bf Proof.}
The proofs of (4.4)-(4.6) are similar as Lemma 3.2. We prove (4.7) as following.

Noting that $J_{1,1}(h)=(h-\frac{1}{6})J_{1,1}^{'}(h)$, we can get $J_{1,1}^{(k)}(h)=0$ $(k\geq2)$, than differentiating (4.5) with respect to $h$ we have
$$
-\frac{1}{2}J_{0,1}^{'}(h)=\frac{3}{2}hJ_{0,1}^{''}(h)+\frac{1}{4}J_{2,1}^{''}(h),
\eqno(4.8)$$
differentiating (4.8) with respect to $h$, we have
$$
hJ_{0,1}^{(3)}(h)=-\frac{4}{3}J_{0,1}^{''}(h)-\frac{1}{6}J_{2,1}^{(3)}(h).
\eqno(4.9)$$
By (4.5), we get that
$$
J_{1,1}(h)=(h-\frac{1}{6})J_{1,1}^{'}(h),~~~J_{1,1}^{''}(h)=0,
$$
taking the derivative of (4.5) twice, we have
$$
\left(\begin{matrix}
     J_{0,1}^{''}(h)\\
     J_{2,1}^{''}(h)
     \end{matrix}\right)
=\left(\begin{matrix}
        &-\frac{3}{4}h~~&-\frac{1}{8}\\
        &\frac{3}{16}h~~&\frac{15}{32}-\frac{3}{2}h
        \end{matrix}\right)
\left(\begin{matrix}
      J_{0,1}^{(3)}(h)\\
      J_{2,1}^{(3)}(h)
      \end{matrix}\right)
\eqno(4.10)$$
Substituting (4.8) and (4.9) into (4.10), we can obtain the (4.7). $\diamondsuit$

Using (4.4) and (4.5), similar to the process of calculation (3.21) and (3.22), we can obtain
$$
\gamma_{1}(h)M^{'}(h)=\gamma_{2}(h)M(h)+F_{1}(h),
\eqno(4.11)$$
where
$$
F_{1}(h)=\alpha_{3}(h)J_{0,0}^{'}(h)+\beta_{3}(h)J_{1,0}^{'}(h)+\eta_{3}(h)J_{0,1}^{'}(h)+\xi_{3}(h)J_{1,1}^{'}(h)+\zeta_{3}(h)J_{2,1}^{'}(h),
$$
with ${\rm deg}\alpha_{3}(h)\leq\left[\frac{n}{3}\right]+\left[\frac{n-2}{3}\right]+1$, ${\rm deg}\beta_{3}(h)\leq\left[\frac{n-1}{3}\right]+\left[\frac{n-2}{3}\right]$, ${\rm deg}\eta_{3}(h)\leq\left[\frac{n-1}{3}\right]+\left[\frac{n-2}{3}\right]+1$, ${\rm deg}\xi_{3}(h)\leq2\left[\frac{n-2}{3}\right]+1$, ${\rm deg}\zeta_{3}(h)\leq\left[\frac{n-3}{3}\right]+\left[\frac{n-2}{3}\right]$. By Lemma 5.1 of [15], we have
$$
\begin{aligned}
\#\{M(h)=0,h\in(0,\frac{1}{6})\}&\leq\#\{\gamma_{1}(h)=0,h\in(0,\frac{1}{6})\}\\
&+\#\{F_{1}(h)=0,h\in(0,\frac{1}{6})\}+1,
\end{aligned}
\eqno(4.12)$$
Let
$$
\Phi_{1}(h)=\alpha_{3}(h)J_{0,0}^{'}(h)+\beta_{3}(h)J_{1,0}^{'}(h),
$$
$$
\Phi_{2}(h)=\eta_{3}(h)J_{0,1}^{'}(h)+\xi_{3}(h)J_{1,1}^{'}(h)+\zeta_{3}(h)J_{2,1}^{'}(h),
$$
similar to the proof of Lemma 3.3, Lemma 3.4 and Lemma 3.5, we have the following Lemmas.

\vskip 0.3 true cm
\noindent{\bf Lemma 4.3.}
If $n\geq3$, for $h\in(0,\frac{1}{6})$, we have
$$
J_{0,0}^{'}(h)\neq0,~~~J_{0,1}^{'}(h)\neq0.
$$
Let
$$
\omega(h)=\frac{J_{1,0}^{'}(h)}{J_{0,0}^{'}(h)},~~\nu(h)=\frac{J_{2,1}^{''}(h)}{J_{0,1}^{'}(h)},~~D_{3}(h)=\frac{3}{64}h(24h-7),
$$
then $\omega(h)$ and $\nu(h)$ satisfy the following Riccati equations respectively
$$
D_{1}(h)\omega^{'}(h)=-k_{12}(h)\omega^{2}(h)+(k_{22}(h)-k_{11}(h))\omega(h)+k_{21}(h),
\eqno(4.13)$$
$$
D_{3}(h)\nu^{'}(h)=-l_{12}(h)\nu^{2}(h)+(l_{22}(h)-l_{11}(h))\nu(h)+l_{21}(h).
\eqno(4.14)$$

\vskip 0.3 true cm
\noindent{\bf Lemma 4.4.}
For $h\in(0,\frac{1}{6})$, there exists polynomials $P_{2}(h)$, $P_{1}(h)$ and $P_{0}(h)$ of $h$ satisfy
${\rm deg}P_{2}(h)\leq m_{2}$, ${\rm deg}P_{1}(h)\leq m_{1}$, ${\rm deg}P_{0}(h)\leq m_{0}$ such that $$L(h)\Phi_{1}(h)=0,$$ where
$m_{2}=\left[\frac{n}{3}\right]+\left[\frac{n-1}{3}\right]+2\left[\frac{n-2}{3}\right]+13$, $m_{1}=m_{2}-1$, $m_{0}=m_{2}-2$,
and
$$
L(h)=P_{2}(h)\frac{d^{2}}{dh^{2}}+P_{1}(h)\frac{d}{dh}+P_{0}(h),
\eqno(4.15)$$

\vskip 0.3 true cm
\noindent{\bf Lemma 4.5.} Suppose that $D_{3}(h)\neq0$ for $h\in(0,\frac{1}{6})$, then we have $L(h)F_{1}(h)=R(h)$, where $L(h)$ is defined as (4.15), and
$$
R(h)=\frac{1}{D_{3}^{2}(h)}\left(Q_{1}(h)J_{0,1}^{'}(h)+Q_{2}(h)J_{2,1}^{''}(h)\right)+Q_{3}(h)J_{1,1}^{'}(h),
\eqno(4.16)$$
$Q_{s}(h)(s=1,2,3)$ are polynomials of $h$ with ${\rm deg}Q_{1}(h)\leq\left[\frac{n}{3}\right]+2\left[\frac{n-1}{3}\right]+3\left[\frac{n-2}{3}\right]+16$, ${\rm deg}Q_{2}(h)\leq\left[\frac{n}{3}\right]+\left[\frac{n-1}{3}\right]+3\left[\frac{n-2}{3}\right]+\left[\frac{n-3}{3}\right]+17$ and ${\rm deg}Q_{3}(h)\leq\left[\frac{n}{3}\right]+\left[\frac{n-1}{3}\right]+4\left[\frac{n-2}{3}\right]+12$.

\vskip 0.3 true cm
\noindent{\bf Lemma 4.6.}
If $n\geq3$, for $h\in(0,\frac{1}{6})$, the number of zeros of $\Phi_{1}(h)$ does not exceed $\left[\frac{n}{3}\right]+2\left[\frac{n-1}{3}\right]+3\left[\frac{n-2}{3}\right]+9$, and the number of zeros of $R(h)$ does not exceed $7\left[\frac{n}{3}\right]+9\left[\frac{n-1}{3}\right]+25\left[\frac{n-2}{3}\right]+\left[\frac{n-3}{3}\right]+105$.

\vskip 0.1 true cm
\noindent{\bf Proof.}
By the same arguments as the proof of Lemma 3.6, we get the estimate of number of $\Phi_{1}(h)$.
Noting
$$
J_{1,1}^{(k)}=0,~~~k\geq2,
$$
let $k_{0}=\left[\frac{n}{3}\right]+\left[\frac{n-1}{3}\right]+4\left[\frac{n-2}{3}\right]+13$,deriving $k_{0}$ times for $R(h)$ with respect to $h$, we have
$$
R^{(k_{0})}(h)=\frac{1}{D_{3}^{k_{0}+2}}\left(\eta_{4}(h)J_{0,1}^{'}(h)+\zeta_{4}(h)J_{2,1}^{''}(h)\right),
\eqno(4.17)$$
where ${\rm deg}\eta_{4}(h)\leq2\left[\frac{n}{3}\right]+3\left[\frac{n-1}{3}\right]+7\left[\frac{n-2}{3}\right]+29$, ${\rm deg}\eta_{4}(h)\leq2\left[\frac{n}{3}\right]+2\left[\frac{n-1}{3}\right]+7\left[\frac{n-2}{3}\right]+\left[\frac{n-3}{3}\right]+30$.
Let
$$
R_{1}(h)=\eta_{4}(h)J_{0,1}^{'}(h)+\zeta_{4}(h)J_{2,1}^{''}(h),
$$
$$
S_{2}(h)=\frac{R_{1}(h)}{J_{0,1}^{'}(h)}=\eta_{4}(h)+\zeta_{4}(h)\nu(h).
$$
Using (4.14), we can end the proof. $\diamondsuit$

\vskip 0.3 true cm
\noindent{\bf Proof for the case of Hamiltonian triangle.}

 Denote $m_{5}=\left[\frac{n}{3}\right]+2\left[\frac{n-1}{3}\right]+3\left[\frac{n-2}{3}\right]+9$, $m_{6}=7\left[\frac{n}{3}\right]+9\left[\frac{n-1}{3}\right]+25\left[\frac{n-2}{3}\right]+\left[\frac{n-3}{3}\right]+105$. For $h\in(0,\frac{1}{6})$, by Lemma 4.4 and Lemma 4.5 we can assume that
$$
P_{2}(\widetilde{h}_{i})=0,~~\Phi_{1}(\overline{h}_{j})=0,~~\widetilde{h}_{i},\overline{h}_{j}\in(-2,0),~~1\leq i\leq m_{2} ,~~1\leq j\leq m_{5}.
$$
Denote $\widetilde{h}_{i}$ and $\overline{h}_{j}$ as $h_{k}^{*}<h_{k+1}^{*}$ for $k=1,2,...,m_{2}+m_{3}$. Let
$$
\Delta_{s}=(h_{s}^{*},h_{s+1}^{*}),~~s=0,1,...,m_{2}+m_{5},
$$
where $h_{0}^{*}=0$ and $h_{m_{2}+m_{5}+1}^{*}=0$. Then $P_{2}(h)\neq0$ and $\Phi_{1}(h)\neq0$ for $h\in\Delta_{s}$ and $L(h)\Phi_{1}(h)=0$. By Lemma 2.3, the solution space of
$$
L(h)=P_{2}(h)\left(\frac{d^{2}}{dh^{2}}+\frac{P_{1}(h)}{P_{2}(h)}\frac{d}{dh}+\frac{P_{0}(h)}{P_{2}(h)}\right)
$$
is a Chebyshev space on $\Delta_{s}$. By Lemma 2.4, $M(h)$ has at most $2+l_{s}$ zeros for $h\in\Delta_{s}$, where $l_{s}$ is the number of zeros of $R(h)$ on $\Delta_{s}$. We obtain
$$
\begin{aligned}
\#\{M(h)=0,h\in\Sigma_{1}\}&\leq \#\{\gamma_{1}(h)=0,h\in(0,\frac{1}{6})\}+\#\{F_{1}(h)=0,h\in(0,\frac{1}{6})\}+1\\
&\leq \#\{R(h)=0,h\in\Sigma_{1}\}+2\cdot the ~number ~of ~the ~intervals ~of ~\Delta_{s}\\
&+the ~number ~of ~the ~end ~points ~of ~\Delta_{s}+\#\{\gamma_{1}(h)=0,h\in(0,\frac{1}{6})\}+1\\
&\leq m_{6}+3(m_{2}+m_{5}+1)+\left[\frac{n-2}{3}\right]+4\\
&\leq 24n+126.
\end{aligned}$$
This ends the proof.

\end{document}